\documentclass[12pt,reqno]{amsart}
\usepackage{amssymb,delarray}
\usepackage{amsfonts}
\usepackage{epsfig}
\usepackage[all]{xy}
\usepackage{amscd}
\usepackage{epigraph}

\makeindex{}

\def\le{\leqslant}
\def\ge{\geqslant}

\def\Kl{\operatorname {Kl}}
\def\Aut{\operatorname {Aut}}

\newtheorem{thm}{Theorem}[section]
\newtheorem{lem}[thm]{Lemma}
\newtheorem{prop}[thm]{Proposition}

\newtheorem{cor}[thm]{Corollary}

{\catcode`\@=11
\gdef\n@te#1#2{\leavevmode\vadjust{%
 {\setbox\z@\hbox to\z@{\strut#1}%
  \setbox\z@\hbox{\raise\dp\strutbox\box\z@}\ht\z@=\z@\dp\z@=\z@%
  #2\box\z@}}}
\gdef\leftnote#1{\n@te{\hss#1\quad}{}}
\gdef\rightnote#1{\n@te{\quad\kern-\leftskip#1\hss}{\moveright\hsize}}
\gdef\?{\FN@\qumark}
\gdef\qumark{\ifx\next"\DN@"##1"{\leftnote{\rm##1}}\else
 \DN@{\leftnote{\rm??}}\fi{\rm??}\next@}}

\begin{document}

\baselineskip=14.pt plus 2pt 

\title[Automorphisms of Galois coverings]{Automorphisms
of Galois coverings of generic $m$-canonical projections}
\author[V.M. Kharlamov and  Vik.S.~Kulikov]{V.M. Kharlamov and  Vik.S. Kulikov}
\address{
UFR de Math\'{e}matiques et IRMA \\ Universit\'{e} Louis Pasteur et
CNRS \\ 7 rue R\'{e}n\'{e}-Descartes
\\ 67084, Strasbourg, Cedex, France}
\email{ kharlam@math.u-strasbg.fr}
\address{Steklov Mathematical Institute\\
Gubkina str., 8\\
119991 Moscow \\
Russia} \email{kulikov@mi.ras.ru}

\thanks{This
research was performed during the stay of the second author in
Strasbourg University under a contract with CNRS and in the
framework of the grant INTAS 05-7805. The first author acknowledges
a support from Agence Nationale de la Recherche through
ANR-05-0053-01, the second author acknowledges a support by RFBR
(05-01-00455), NWO-RFBR 047.011.2004.026
(RFBR 05-02-89000-$\text{NWO}_{-}a$), RFBR
06-01-72017-$\text{MNTI}_{-}a$, and RUM1-2692-MO-05.} \keywords{}

\begin{abstract}
The automorphism group of the Galois covering induced by a
pluri-canonical generic covering of a projective space is
investigated. It is shown that by means of such coverings one
obtains, in dimensions one and two, serieses of specific actions of
the symmetric groups $S_d$ on curves and surfaces not deformable to
an action of $S_d$ which is not the full automorphism group. As an
application, new DIF$\ne$DEF examples for $G$-varieties in complex
and real geometry are given.
\end{abstract}

\maketitle

\epigraph{\scriptsize Perhaps the simplest combinatorial entity is
the group of the $n!$ permutations of $n$ things. This group has a
different constitution for each individual number $n$. The question
is whether there are nevertheless some asymptotic uniformities
prevailing for large $n$ or for some distinctive class of large $n$.
Mathematics has still little to tell about such problems.} {\small
\it H.~Weyl, Philosophy of Mathematics and Natural Science, Appendix
D, Princeton Univ. Press, 1949.}

\setcounter{tocdepth}{2}


\def\st{{\sf st}}

\section*{Introduction.}

\subsection{Terminology conventions}
By a {\it covering}\, we understand a branched covering, that is a
finite morphism $f:X\to Y$ from a normal projective variety $X$ onto
a non-singular projective variety $Y$, all being defined over the
field of complex numbers $\mathbb C$.  To each covering $f$ we
associate the branch locus $B\subset Y$, the ramification locus
$R\subset f^{-1}(B)\subset X$, and the unramified part $X\setminus
f^{-1}(B)\to Y\setminus B$ (which is the maximal unramified
subcovering). As is usual for unramified coverings, there is a
homomorphism $\psi$ which acts from the fundamental group $\pi_1(
Y\setminus B)$ to the symmetric group $S_d$ on $d$ elements, where
$d$ is the degree of $f$. This homomorphism (called monodromy of
$f$) is defined by $f$ uniquely up to inner automorphisms of $S_d$;
reciprocally, according to Grauert-Remmert-Riemann extension theorem
the conjugacy class of $\psi$ defines $f$ up to an isomorphism. The
image $G\subset S_d$ of $\psi$ is a transitive subgroup of $S_d$.

As is well known, for each covering $f$ there is a unique, up to
isomorphism, minimal Galois covering $\widetilde f:\widetilde X\to
Y$ which is factorized through $f$, $\widetilde f=f\circ h$, by
means of a Galois covering $h:\widetilde X\to X$. The covering
$\widetilde f$ is called the {\it Galois expansion of $f$.} The
characteristic, minimality, property of the Galois expansion
$\widetilde f$ is that any Galois covering which is factorized
through $f$ can be factorized through $\widetilde f$. The  Galois
expansion of $f$ can be obtained by Grauert-Remmert-Riemann
extension theorem from the non-diagonal component of the fibred
product over $Y$ of $d$ copies of the unbranched part $X\setminus
f^{-1}(B)\to Y \setminus B$ of $f$. In particular, the Galois
group $Gal(\widetilde X/Y)$ of the covering $\widetilde
f:\widetilde X\to Y$ is naturally identified with
$G=\psi(\pi_1(Y\setminus B))$ as above.

In the present article, we study actions of finite groups on the
Galois expansions of generic coverings of the projective spaces. To
give a precise definition of a generic covering we need to introduce
first some preliminary definitions concerning the actions of
symmetric groups.

Let $I$ be a finite set consisting of $|I|=d$ elements and let
$I_1\cup \dots \cup I_k=I$ be a partition of $I$, $|I_i|=d_i\geq 1$,
$\sum_{i=1}^k d_i=d$. Such a partition defines a unique, up to
conjugation, imbedding of $S_{d_1}\times \dots\times S_{d_k}$ in
$S_d$ which we call a {\it standard imbedding}. A representation
$S_{d_i}\subset GL(V_i)$, where $V_i$ is a vector space over
$\mathbb C$, will be called a {\it standard representation of rank}
$d_i-1$ if there is a base $e_1,\dots, e_{r_i}$ of $V_i$, $r_i=\dim
V_i\geq d_i-1$, such that the action $\sigma_{(j,j+1)}$ of a
transposition $(j,j+1)
\in S_{d_i}$ is given by
$$\sigma_{(j,j+1)}(e_l)=\left\{
\begin{array}{ll}
 e_l  & \, \text{if} \,\, l\neq j,j+1, \\
 e_{j+1} & \, \text{if} \,\, l=j,
\end{array} \right. $$
for $j\neq d_i-1$, and by
$$\sigma_{(d_i-1,d_i)}(e_l)=\left\{
\begin{array}{cl}  e_l
 & \, \text{if} \,\, l\neq d_i-1, \\
-\sum_{s=1}^{d_i-1 } e_{s} & \, \text{if} \,\, l=d_i-1,
\end{array} \right. $$
in the remaining case. A collection of standard representations of
symmetric groups $S_{d_i}\subset GL(V_i)$ of ranks $d_i-1$,
$i=1,\dots, k$, defines a representation of $S_{d_1}\times
\dots\times S_{d_k}\subset GL(V)$ with $V=V_{i_1}\oplus\dots\oplus
V_{i_k}$ which we call  a {\it standard representation of the
product $S_{d_1}\times \dots\times S_{d_k}$ of rank} $\sum d_i-k$.
As is easy to see, if $S_{d_1}\times \dots\times S_{d_k}\subset
GL(V)$ is a  standard representation of rank $\sum d_i-k$, then the
codimension in $V$ of the subspace consisting of the vectors fixed
under the action of $S_{d_1}\times \dots\times S_{d_k}$ is equal to
$\sum d_i-k$.

Let the group $S_d$ act on a smooth projective manifold $Y$. We say
that the action of $S_d$ on $Y$ is {\it generic} if the stabilizer
$St_a\subset S_d$ of each point $a\in Y$ is a standard imbedding of
a product of symmetric groups and the action induced by $St_a$ on
the tangent space $T_aY$ is a standard representation (both the
product and the representation depending on $a$). According to this
definition, if the action of $S_d$ on $Y$ is generic, then the
factor-space $Y/S_d$ is a smooth projective manifold.

A covering $f:X\to \mathbb P^{\dim X}$ of degree $d$ is called {\it
generic} if the Galois group $G=Gal(\widetilde X/\mathbb P^{\dim
X})$ of the Galois expansion of $f$ is the full symmetric group
$S_d$, the varieties $X$ and $\widetilde X$ are smooth, and the
action of $G$ on $\widetilde X$ is generic. From this definition it
follows that, for any generic covering $f$ of degree $d$, the group
$Gal(\widetilde X/\mathbb P^{\dim X})$ is the full symmetric group
$S_d$ and the subgroup $Gal(\widetilde X/X)$ coincides with
$S_{d-1}\subset S_d$.

If $X$ is non-singular and $\dim X=1$, then a covering $f:X\to
\mathbb P^1$ branched over $B\subset \mathbb P^1$ is generic if and
only if $|f^{-1}(b)|=\deg f-1$ for any $b\in B$. Furthermore, in the
case of generic covering at each point $\widetilde b\in \widetilde
f^{-1}(b), b\in B$, the stabilizer group $St_{\widetilde b} \subset
S_d=Gal(\widetilde X/\mathbb P^1)$ is generated by a transposition,
$St_{\widetilde b}= S_2$.

If $X$ is non-singular and $\dim X=2$, then a covering $f:X\to
\mathbb P^2$ is generic if and only if the following conditions are
satisfied: $f$ is branched over a cuspidal curve $B\subset \mathbb
P^2$; $|f^{-1}(b)|=\deg f-1$ for any nonsingular point $b$ of $B$,
and $|f^{-1}(b)|=\deg f-2$ if $b$ is a node or a cusp of $B$. In the
case of generic covering at each point $\widetilde b\in \widetilde
f^{-1}(b)$, $b\in B$, the stabilizer group $St_{\widetilde b}\subset
S_d=Gal(\widetilde X/\mathbb P^2)$ is generated: by a transposition
if $b$ is a nonsingular point of $B$, and then $St_{\widetilde
b}=S_2$; by two non-commuting transpositions if $b$ is a cusp of
$B$, and then $St_{\widetilde b}=S_3$; by two commuting
transpositions if $b$ is a node of $B$, and then $St_{\widetilde b}=
S_2\times S_2$ (for a detailed exposition see Subsection
\ref{sec31}).

Whatever is the dimension of a generic covering, the automorphism
group $Aut(\widetilde X)$ of the manifold $\widetilde X$ contains
the symmetric group $S_d$, but as the following examples show, one
can not expect that $Aut(\widetilde X)$ and $S_d$ will necessarily
coincide.

As a first example, let us pick a generic covering $f_1:Y = \mathbb
P^1\to\mathbb P^1$ of degree $d+1$ and denote by $\widetilde
f:\widetilde Y \to \mathbb P^1$ the Galois expansion of $f_1$,
$\widetilde f=f_1\circ h_1$. Then $Gal(\widetilde Y/ \mathbb
P^1)=S_{d+1}$, $Y=\widetilde Y/S_{d}$, and $h_1:\widetilde Y \to
Y=\mathbb P^1$ is a Galois covering with $Gal(\widetilde
Y/Y)=S_{d}$. The covering $h_1 : \widetilde Y \to Y=\mathbb P^1$ can
be considered as the Galois expansion $\widetilde X=\widetilde Y\to
Y=\mathbb P^1$ (with Galois group $Gal(\widetilde X/ \mathbb
P^1)=S_{d}$) of a covering $f:X\to Y=\mathbb P^1$, $X=\widetilde
X/S_{d-1}$. The latter is a generic covering of degree $d$, and,
now, if we start from $f:X\to \mathbb P^1$, we obtain that its
Galois group $Gal(\widetilde X/\mathbb P^1)=S_d$ does not coincide
with $Aut(\widetilde X)$, since $Aut(\widetilde X)=Aut(\widetilde
Y)$ contains at least the group $S_{d+1}$.

Another example can be obtained as follows. Let $f_1:Y=\mathbb
P^1\to\mathbb P^1$ be a generic degree $d$ covering branched over
$B_1\subset \mathbb P^1$. Let us choose two points $x,y\in  \mathbb
P^1$ not belonging to $B_1$, and let $f_2:Z=\mathbb P^1\to\mathbb
P^1$ be a cyclic covering of degree $p$ branched at $x$ and $y$.
Consider the fibred product $X=Y\times_{\mathbb P^1}Z$ and its
projection $f:X\to Z= \mathbb P^1$ to the second factor. It is easy
to see that $f$ is a generic covering and $Aut(\widetilde X)$
contains $Gal(\widetilde X/\mathbb P^1)\times \mathbb Z/ p\mathbb
Z$.

\subsection{Principal results}
The aim of our research is to give numerical conditions for a
generic covering $f: X \to \mathbb P^{\dim X}$ which ensure that
$Aut(\widetilde X)= Gal(\widetilde X/\mathbb P^{\dim X})$ and which
are preserved under any deformation of the Galois expansion.

To state the results obtained we need to introduce one more
auxiliary notion: a covering  $f:X\to \mathbb P^{\dim X}$ is said to
be (numerically) {\it $m$-canonical} if it is given by $\dim X+1$
sections of a line bundle numerically equivalent to the $m$-th power
$K_X^{\otimes m}$ of the canonical bundle $K_X$ of $X$ and these
sections have no common zeros.

Certainly, $m$-canonical coverings exist only if the Kodaira
dimension of $X$ coincides with its dimension. If $\dim X=2$, then,
in addition, $X$ should be minimal and it should not contain any
$(-2)$-curve. If $\dim X=1$, then its genus should be greater or
equal to $2$. Reciprocally, as is well known, any curve of genus
$g\ge 2$ possesses a $m$-canonical covering for  $m\ge 1$, and, as
is shown in \cite{Ku-Ku}, any minimal surface of general type
containing no $(-2)$-curves also possesses a $m$-canonical covering
at least for $m\ge 10$.

\begin{thm} \label{main1} Let $X$ be a curve of genus $g\ge 2$
and $\widetilde f:\widetilde X\to \mathbb P^1$ be the Galois
expansion of a $m$-canonical generic covering $f: X\to \mathbb P^1$.
If $m(g-1) \geq 500$, the Galois group $Gal(\widetilde X/\mathbb
P^1)$ is the full automorphism group of $\widetilde X$.
\end{thm}

\begin{thm} \label{main2} Let $X$ be a surface of general type,
and assume that it possesses a $m$-canonical generic covering $f:
X\to \mathbb P^2$, $m\geq 2$. If $m^2K_X^2\geq 2 \cdot 84^2$ and
$\widetilde f:\widetilde X\to \mathbb P^2$ is the Galois expansion
of
$f$, the Galois group $Gal(\widetilde X/\mathbb P^2)$ is the full
automorphism group of $\widetilde X$.
\end{thm}

As a consequence,\ the $G$-curves like in Theorem \ref{main1} and
the $G$-surfaces like in Theorem \ref{main2} provide infinitely many
examples of {\it saturated} connected components in the moduli space
of $G$-varieties, $G = S_d$, where a component is called saturated
if for any $G$-variety representing a point of this component $G$ is
the full automorphism group of $V$. (It may be worth pointing some
easy series of saturated components with $G\ne S_d$, namely, the
components given by curves and surfaces with the automorphism groups
of maximal order, that is $84(g-1)$ for curves, and $42^2K^2$ for
surfaces. One can mention also deformation rigid varieties with
nontrivial automorphism group.)

As another application of the above theorems, we give
counter-examples to a Dif$=$Def problem for complex and real
$G$-varieties. Namely, we construct pairs of complex (respectively,
real) varieties $V_1, V_2$ such that the actions of $\Aut V_1$ and
$\Aut V_2$ (respectively, $\Kl V_1$ and $\Kl V_2$; here, $\Kl V$ is
the group formed by the regular isomorphisms $X\to X$ and $X\to\bar
X$) are diffeomorphic but not deformation equivalent. Up to our
knowledge, such examples, specially at the real setting, are new.
(It may be worth noticing, that in \cite{Kh-Ku1} in our
counter-examples  to the Dif$=$Def problem for real structures the
surfaces have diffeomorphic real structures, but the actions of the
Klein group on these surfaces are not diffeomrophic.)

\subsection{Contents of the paper}
The proof of theorems \ref{main1}
and \ref{main2} consists of two parts. In the beginning (Section 1),
by methods of group theory, we investigate minimal expansions of the
symmetric groups to restrict the number of possible cases, and then
the possible cases are investigated by geometric methods (Section 2
for Theorem \ref{main1} and Section 3 for Theorem \ref{main2}).
Section 4 contains the applications mentioned above.

\section{Minimal expansions of symmetric groups.}
\subsection{Preliminary
definitions} To formulate group theoretic statements which we use in
the proof of Theorems \ref{main1} and \ref{main2}, we need to
introduce few preliminary definitions. We say that a group $G$
containing the symmetric group $S_d$ satisfies the {\it minimality
property} if there is no any proper subgroup $G_1$ of $G$ which
contains $S_d$ and does not coincide with $S_d$, and call such a
group $G$ a {\it minimal expansion of $S_d$.}

An imbedding $\alpha : S_d\subset S_{d+2}$ is called {\it
quasi-standard} if the image $\alpha (\sigma_{i,j})$ of each
transposition $\sigma_{i,j} =(i,j)\in S_d$, $1\leq i,j\leq d$, is
the product $\alpha (\sigma_{i,j})=(i,j)(d+1,d+2)$ of two
transpositions $(i,j)$ and $(d+1,d+2)$ in $S_{d+2}$. Note that for
the quasi-standard imbedding the image of $S_d$ is contained in the
alternating subgroup $A_{d+2}$ of $S_{d+2}$. This imbedding, $\alpha
: S_d\hookrightarrow A_{d+2}$, is called {\it standard}.

\begin{prop} \label{prop1}
Let $G$ be a minimal expansion of the symmetric group $S_d$ of index
$k=(G:S_d)$. Assume that $k \leq cd^n$, where either {\rm (i)}
$c=63$ and $n=1$, or {\rm (ii)} $c=(4\cdot 42)^2$ and $n=2$.
\par
If $d\geq
\max (2c, 1000)$, then $G$ is one of the following groups:
\begin{itemize}
\item[($1$)] $G=S_d\times \mathbb Z/p\mathbb Z$, $p$ is a prime number,
$p\leq cd^n$;
\item[($2$)] $G=A_d\rtimes D_r$, where
$3\leq r\leq cd^n$, $r$ is odd, $D_r$ is the dihedral group given by
presentation
$$
D_r=\langle \sigma, \tau \, |\, \sigma^2=\tau^2=(\sigma\tau)^r=1\rangle,
$$
and the action {\rm (}by conjugation{\rm )} of $\sigma$ and $\tau$
on $A_d$ coincides with the action of the transposition $(1,2)\in
S_d$ on $A_d\subset S_d$; \item[($3$)] $G=S_{d+1}$ is the
symmetric group; \item[($4$)] $G=A_{d+2}$ is the alternating
group, the imbedding of $S_d$ in $A_{d+2}$ is a standard one, and
such expansions can appear only under assumption {\rm (ii)}.
\end{itemize}
\end{prop}

The rest of this section is devoted to the proof of Proposition \ref{prop1}.

\proof A priori, one of the following two cases occurs.

\textsf{Case $I$.} {\it The group $G$ contains a non-trivial normal
subgroup.}

\textsf{Case $II$.} {\it The group $G$ is simple.}

Since $d>6$, the group $S_d$ has the unique non-trivial normal
subgroup, namely, the alternating group $A_d$ and, consequently,
Case $I$ can be subdivided into the following subcases, where $N$
denotes a nontrivial normal subgroup of $G$.

\textsf{Case $I_1$.} $S_d\subset N$.

\textsf{Case $I_2$.} $N\cap S_d=\{ 1\}$.

\textsf{Case $I_3$.} $N\cap S_d=A_d$.

In its turn, Case $I_3$ can be subdivided into two subcases.

\textsf{Case $I_{31}$.} {\it $A_d$ is a normal subgroup of $G$}.

\textsf{Case $I_{32}$.} {\it $A_d$ is not a normal subgroup of $G$}.

\subsection{Analysis of Case $I_1$}\label{I1}
It follows from the minimality property that $S_d=N$. Let $g_1$ be
an arbitrary element of $G\setminus S_d$. The conjugation by $g_1$
induces an automorphism of $S_d$. Since $d\geq 7$, any automorphism
of $S_d$ is inner. Therefore, there is $g_2\in S_d$ such that
$g=g_1g_2$ commutes with all elements of $S_d$. Hence, once more by
the minimality property, the group $G$ splits into the direct
product of $S_d$ and the cyclic group $\langle g\rangle$ generated
by $g$. Moreover, the order of $g$ is a prime number $p$.

\subsection{Analysis of Case $I_{31}$}\label{I31}
According to subsection \ref{I1} we can assume that $S_d$ is not a
normal subgroup of $G$. Therefore, there is $g\in G$ such that
$S^{\prime}_d=g^{-1}S_dg$ does not coincide with $S_d$ (but, it is
isomorphic to $S_d$).

Since the group $A_d$ is a normal subgroup of $G$, we have
$A_d\subset S^{\prime}_d\cap S_d$. Furthermore, for any
transposition $\sigma\in S_d$ the element $\tau=g^{-1}\sigma g$
(which we call a transposition in $S^{\prime}_d$) does not belong
to $S_d$. Thus, by the minimality property, it follows that the
group $G$ is generated by the elements of $A_d$ and any two
transpositions $\sigma\in S_d$ and $\tau\in S^{\prime}_d$.
Moreover, since conjugating by elements of $S_d$ (respectively,
$S_d^\prime$) provides the full automorphism group $Aut(A_d)$ of
$A_d$, we can choose the two generating transpositions $\sigma\in
S_d$ and $\tau\in S^{\prime}_d$ in a way that $\sigma\tau$
commutes with all elements of $A_d$. And above all, we can assume
that the action (by conjugation) of $\sigma$ and $\tau$ on $A_d$
coincides with the action of the transposition $(1,2)\in S_d$.

Denote by $H$ a subgroup of $G$ generated by  $\sigma$ and $\tau$.
Then $H$ is isomorphic to a dihedral group
$$D_r=\langle \sigma, \tau\, |\, \sigma^2=\tau^2=(\sigma\tau)^r=1\rangle$$
for some $r\in\mathbb N$.

As is known, any element $g\in D_r$ either belongs to the cyclic
subgroup generated by  $\sigma\tau$ or is conjugate to $\sigma$ or
$\tau$. Therefore $A_d\cap H$ is a subgroup of $\langle
\sigma\tau\rangle$, and since the element $\sigma\tau$ commutes with
all elements of $A_d$ and $A_d$ has trivial center, we conclude that
$A_d\cap H=\{ 1\}$. In addition, $A_d$ is a normal subgroup of $G$
and $G$ satisfies the minimality property, which implies
$$
G=A_d\rtimes H\simeq A_d\rtimes D_r.
$$
Moreover, $r$ is odd, since $\sigma$ and $\tau$ are conjugate in $G$
and, hence, in $D_r$.

\subsection{Analysis of Case $I_2$}
It follows from the minimality property that in this case the group
$G$ is isomorphic to a semi-direct product $N \rtimes S_d$.

If $N$ is not a simple group, then we can find a minimal non-trivial
normal subgroup $N_1$ of $N$. Note, first of all, that $N_1$ can not
be a normal subgroup of $G$, since $G$ satisfies the minimality
property. Therefore, the set of subgroups of $G$ conjugated to $N_1$
contains more than one element. Let $\{N_1, \dots, N_s\}$ be the set
of subgroups conjugate to $N_1$ in $G$, $s\ge 2$. Each $N_i, 1\le
i\le s,$ is contained in $N$, since $N$ is a normal subgroup of $G$.
Moreover, each of them is a normal subgroup of $N$, since $N_1$ is
normal in $N$ and conjugation by any element of $G$ induces an
isomorphism of $N$. Besides, the action of $S_d$ on the set $\{N_1,
\dots, N_s\}$ is transitive and this set is an orbit of the action
of $S_d$ by conjugation on the whole set of subgroups of $G$.

Let us show that $N\simeq  N_1\times\dots\times N_{s_1}$ for some
$s_1\leq s$ (maybe after a renumbering the groups $N_i$). Note,
first, that $N_i\cap N_j=\{ 1\}$ for $i\neq j$. Indeed,
$N_1,\dots,N_s$ are minimal normal subgroups of $N$ and the
intersection $N_i\cap N_j$, as an intersection of any two normal
subgroups, is a normal subgroup. Therefore, $[N_1,N_2]\subset
N_1\cap N_2=\{ 1\}$ so that the subgroup $N_1N_2$ generated in $N$
by the elements of the groups $N_1$ and $N_2$ is isomorphic to
$N_1\times N_2$, and this subgroup $N_1N_2$ is also a normal
subgroup of $N$. By induction, assume that for some $i<s$ the
subgroup $N_{1,i}=N_1 \dots N_{i}$ of $N$ is normal and isomorphic
to $N_1\times\dots\times N_{i}$. Then, either $N_{1,i}\cap
N_{i+1}=\{ 1\}$, or $N_{i+1}\subset N_{1,i}$ (here, once more, we
use the observation that $N_{i+1}$ is a minimal normal subgroup of
$N$). If $N_{1,i}\cap N_{i+1}={1}$, then $$N_{1,i+1}=N_1\cdot .\,
.\, . \cdot N_{i+1}\simeq N_1\times\dots\times N_{i+1},$$ since
$[N_{1,i}, N_{i+1}]\subset N_{1,i}\cap N_{i+1}={1}$. This
inductive procedure stops at some subgroup $N_{1,s_1}\subset N$
which, being normal in $N$ and invariant under the action of $S_d$
by conjugation, is a normal subgroup of $G$. Therefore,
$N_{1,s_1}$ coincides with $N$. As a consequence, $N\simeq
N_1\times\dots\times N_{s_1}$.

The groups $N_1,\dots, N_s$, which are isomorphic to each other, are
simple. Indeed, if $N_1$ is not simple, there exists a non-trivial
normal subgroup $\widetilde N_1$ of $N_1$, so that the group
$\widetilde N_1\times \{ 1\}\times\dots\times \{ 1\}$ is a normal
subgroup of $N$, which is impossible, since, by assumption, $N_1$ is
a minimal normal subgroup of $N$.

Next, let us show that $s_1=s$ if $N_1$ is a non-abelian simple
group. Here, we are reasoning by contradiction and suppose that
$s_1<s$. Then, $N_{s_1+1}\subset  N_1\times\dots\times N_{s_1}$, the
projection of $N_{s_1+1}$ to each of the factors is either an
isomorphism or the trivial homomorphism, and at least two
projections are isomorphisms. Without loss of generality, we can
assume that the first two projections are isomorphisms. So, if an
element $g\in N_{s_1+1}$ is written as a product $g=g_1g_2\dots
g_{s_1}$ of elements $g_i\in N_i$, then $g_1$ is uniquely determined
by $g_2$. On the other hand, $N_{s_1+1}$ is a normal subgroup of
$N$, which implies that for any $g=g_1g_2\dots g_{s_1}\in N_{s_1+1}$
and any $h\in N_1$ the products $h^{-1}gh=(h^{-1}g_1h)g_2\dots
g_{s_1}$ belong to $N_{s_1+1}$. Contradiction.

If  $N_1$ is a non-abelian simple group, then conjugation by the
elements of $S_d$ permutes $N_1,\dots, N_s$ transitively and
defines a homomorphism $\psi : S_d\to S_s$. Since $d\geq 7$, there
are only two possibilities: either $s\leq 2$ and $A_d\subset
\ker\psi$, or $\psi$ is an imbedding and, therefore, $s\geq d$.

Suppose that $N_1$ is a non-abelian simple group and $s= 2$. Then,
conjugation by the elements of $A_d$ defines a homomorphism $\psi
: A_d\to Aut(N_1)$.
Since the outer automorphism group $Out(N_1)=Aut(N_1)/Inn (N_1)$ is
solvable for a simple group $N_1$ (see \cite{G}, Theorem 4.240) and
the group $A_d$ is simple, we see that $\psi(A_d)\subset Inn(N_1)$
and either $\psi(A_d)=1$, or $\psi : A_d\to Inn(N_1)$ is an
imbedding into the inner automorphism group. If $\psi(A_d)=1$ , any
element of $N_1$ commutes with the elements of $A_d$, and therefore
$A_d$ is a normal subgroup of $G$, since $G$ is generated by the
elements of $S_d$ and an element $h\in N_1$, $h\neq 1$, and since
$h$ commutes with the elements of $A_d$. But such a case is already
considered above (Case $I_{31}$). If $\psi : A_d\to Inn(N_1)$ is an
imbedding, then the order of $N_1$ is greater than $\frac{1}{2}d!$
and, consequently,
$$k=(G:S_d)\geq (\frac{1}{2}d!)^2,$$ which is impossible since, by
assumption, $d\geq 1000$ and $k$ is not greater than either $63d$
or $(4\cdot 42)^2d^2$.

Suppose, next, that $N_1$ is a non-abelian simple group and $s\geq
d$. Then (as is well known, $A_5$ is the smallest non-abelian
simple group)
$$k=(G:S_d) =|N|^s\geq |N|^d\geq 60^d$$
which is impossible by the same reason as above.

Now, let $N_1$ be an abelian simple group. Then $N\simeq
N_1\times\dots\times N_{s_1}$ and again conjugation by the
elements of $A_d$ defines a homomorphism $\psi : A_d\to Aut(N_1)$.
As above, if  $\psi(A_d)=1$, then $A_d$ is a normal subgroup of
$G$ and this case is already considered.

If $\psi : A_d\to Aut(N_1)$ be an imbedding, then Lemma \ref{dim2}
(see below) implies that $s_1\geq [\frac{d}{4}]$. Therefore,
$$k=(G:S_d) =|N_1|^{s_1}\geq 2^{[\frac{d}{4}]},$$ which
contradicts the assumption that $d\geq 1000$ and $k$ is not
greater than either $63d$ or $(4\cdot 42)^2d^2$.

The rest of the subsection is devoted to a proof of Lemma
\ref{dim2}, which is based, in its turn, on the following Lemma.

\begin{lem} \label{dim}
Let $F$ be a finite field of characteristic $p$ and let
$H$ be a subgroup of $PGL(F,n)$ isomorphic to the alternating group
$A_{4d_1}$, $d_1\in \mathbb N$. Then $n\geq d_1$.
\end{lem}

\proof If $d_1\leq 2$ the statement is obvious. Assume that it  is
true for $d_1\leq k$ and put $d_1=k+1$.

Consider the natural epimorphism $\psi:GL(F,n)\to PGL(F,n)$. The
kernel of $\psi$ consists of scalar matrices $\lambda \text{Id}$,
$\lambda \in F^*$. Therefore the group $\psi^{-1}(A_{4(k+1)})$ is a
central extension of the group $A_{4(k+1)}$.

Denote by $\Xi$ the set of 4-tuples $\{i_1,i_2,i_3,i_4\}$ of
pairwise distinct integers $1\leq i_j\leq 4(k+1)$ and for each
$I\in\Xi$ denote by $x_I$ the permutation $(i_1,i_2)(i_3,i_4)\in
A_{4(k+1)}$. These permutations $x_I$, $I\in\Xi$, generate the group
$A_{4(k+1)}$ and, since $4(k+1)>8$, each two of them are conjugate
in $A_{4(k+1)}$. For each pair $I_1,I_2\in\Xi$ let us choose a word
$w_{I_1,I_2}$ in $x_I$ such that
$x_{I_1}=w_{I_1,I_2}^{-1}x_{I_2}w_{I_1,I_2}$ in $A_{4(k+1)}$, then
pick elements $\hat x_I\in \psi^{-1}(x_{I})\subset
\psi^{-1}(A_{k+1})$ and put
$$\widetilde x_{I_0}=
\hat x_{I_0} \, \, \text{for}\, \, I_0=\{4k+1,4k+2,4k+3,4k+4\},
$$
$$
\widetilde x_{I}=
\hat w_{I,I_0}^{-1}\widetilde x_{I_0}\hat w_{I,I_0}\, \,
\text{for}\, \, I\ne I_0,
$$
where $\hat w_{I,I_0}$ is obtained from the word $w_{I,I_0}$ by
substitution the elements $\hat x_I$ instead of $x_I$.

For each $I\in\Xi$ we have $\widetilde x_I=\mu _I \hat x_I$, where
$\mu_I\in \ker \psi$ are elements of the center of $GL(F,n)$.
Therefore,
$$\widetilde x_{I}=
\widetilde w_{I,I_0}^{-1}\widetilde x_{I_0}\widetilde w_{I,I_0},$$
where $\widetilde w_{I,I_0}$ is obtained from the word $w_{I,I_0}$
by substitution the elements $\widetilde x_I$ instead of $x_I$. On
the other hand, $x_I^2=1$ for each $I\in\Xi$ which implies that
$\widetilde x_I^2=\lambda_I \in \ker \psi$. Since $\widetilde x_I$
are conjugate to each other and $\widetilde x_I^2=\lambda_I$ are the
elements of the center, all $\lambda_I$ should be equal to each
other. We denote this element by $\lambda$.

Consider the group $GL(F,n)$ as
a subgroup of $GL(\overline F,n)$, where
$\overline F$ is the algebraic closure of the field $F$, and denote by
$\widetilde A_{4(k+1)}$ the subgroup of $GL(\overline F,n)$ generated by
$\widetilde x_I, I\in \Xi$, and the elements belonging to the center of
$GL(\overline F,n)$. It is easy to see that $\widetilde A_{4(k+1)}$ is
a central extension of the group $A_{4(k+1)}$.

In $\widetilde A_{4(k+1)}$ there is an element $\mu$ such that
$\mu^2=\lambda^{-1}$. Put $y_I=\mu \widetilde x_I$. Then $y_I^2=1$
and all the $y_I$ are conjugate to each other:
$y_{I}=v_{I,I_0}^{-1}y_{I_0}v_{I,I_0}$, where $v_{I,I_0}$ is
obtained from the word $w_{I,I_0}$ by substitution the elements
$y_I$ instead of $x_I$. Furthermore, for each $I=\{
i_1,i_2,i_3,i_4\}\in\Xi$ with   $1\leq i_j\leq 4k$, we have
$x_Ix_{I_0}=x_{I_0}x_I$ which implies that
$y_Iy_{I_0}=\mu_Iy_{I_0}y_I$ for some $\mu_I$ belonging to the
center of $GL(\overline F,n)$. Since $y_I^2=Id$, then $\mu_I=\pm
Id$, and since all the $y_I$ are conjugate to each other by words
depending on the elements $y_I$, all $\mu_I$ should be equal to each
other. As a result, all $\mu_I$ with $I\in\Xi$ are equal to either
$\mu=Id$ or $\mu= -Id$.

Let us show that $\mu=Id$. Consider the elements $y_{1,2,3,4}$,
$y_{1,2,5,6}$, and put $\widetilde
y_{3,4,5,6}=y_{1,2,3,4}y_{1,2,5,6}$. We have $\widetilde
y_{3,4,5,6}=\lambda  y_{3,4,5,6}$, where $\lambda$ is a central
element, since  $x_{3,4,5,6}=x_{1,2,3,4}x_{1,2,5,6}$. Also, we
have $y_{I_0}y_{1,2,3,4}=\mu y_{1,2,3,4}y_{I_0}$ and
$y_{I_0}y_{1,2,5,6}=\mu y_{1,2,5,6}y_{I_0}$. Hence, on the one
hand,
$$y_{I_0}\widetilde y_{3,4,5,6}=y_{I_0}\lambda y_{3,4,5,6}=
\lambda \mu y_{3,4,5,6}y_{I_0}=\mu \widetilde y_{3,4,5,6}y_{I_0}$$
and, on the other hand,
$$y_{I_0}\widetilde y_{3,4,5,6}=y_{I_0}y_{1,2,3,4}y_{1,2,5,6}=
\mu^2  y_{1,2,3,4}y_{1,2,5,6}y_{I_0}=\mu^2 \widetilde y_{3,4,5,6}y_{I_0}.$$
Therefore, $\mu=Id$.

Denote by $\overline A_{4(k+1)}$ the subgroup of $GL(\overline F,n)$
generated by $y_I, I\in\Xi$. Obviously, its image in $PGL(\overline
F,n)$ is $A_{4(k+1)}$. Consider a subgroup $\overline A_{4k}$ of
$\overline A_{4(k+1)}$ generated by the elements $y_I$, $I=\{
i_1,i_2,i_3,i_4\}\in\Xi$ with   $1\leq i_j\leq 4k$. The elements of
$\overline A_{4k}$ commute with $y_{I_0}$, and the image of
$\overline A_{4k}$  in $PGL(\overline F,n)$ is $A_{4k}$.

Assume, first, that the characteristic $p\neq 2$. Then the vector
space $V=\overline F^n$ splits into a direct sum $E_+\oplus E_{-}$
of two eigen-spaces corresponding to the eigen-values $\pm 1$ of
$y_{I_0}$. Since $y_{I_0}$ does not belong to the center, $\dim
E_+\geq 1$ and $\dim E_{-}\geq 1$. Since the elements of $\overline
A_{4k}$ and $y_{I_0}$ commute, the both eigen-spaces $E_{\pm}$ are
invariant under the action of $A_{4k}$. This action is non-trivial
on at least one of these subspaces, say, on $E_+$. Furthermore,
since $A_{4k}$ is a simple group, this action induces an imbedding
of $A_{4k}$ into $PGL(E_+)$. Therefore, $\dim E_+\geq k$ and, as a
result, $\dim V\geq k+1=d_1$.

Suppose now that $p=2$. Then, the subspace $E= \{ v\in V \, |\,
y_{I_0}(v)=v\}$ of $V$ is invariant  under the action of $A_{4k}$
and it is of dimension $\dim E<\dim V$. If the action of $A_{4k}$ on
$E$ is non-trivial, then $n=\dim V> \dim E \geq k$, that is, $n\geq
k+1=d_1$.

To end the proof, let us show that if the action of $A_{4k}$ on $E$
is trivial, then the induced action of $A_{4k}$ on $V/E$ is
non-trivial. Indeed, if the both actions are trivial, then  we can
choose a basis in $V$ such that in this basis each $y\in A_{4k}$ can
be represented by a matrix of shape
$$
y=
\left( \begin{array}{cc}
Id_a & A \\
0 & Id_b \end{array}
\right) ,
$$
where $a=\dim E$, $b=\dim V- a$, $A$ is a $(a\times b)$-matrix,
and $0$ is the zero $(b\times a)$-matrix. But, it is impossible,
since the such matrices form an abelian group, while the group
$A_{4k}$ is non-abelian. Thus, we conclude that the action of
$A_{4k}$ on $V/E$ is non-trivial, and, hence, $n=\dim V>\dim
V/E\geq k$, that is,  $n\geq k+1=d_1$. \qed

\begin{lem} \label{dim2}
Let $F$ be a finite field of characteristic $p$ and let
$H$ be a subgroup of $GL(F,n)$ isomorphic to the alternating group
$A_{4d_1}$, $d_1\in \mathbb N$. Then $n\geq d_1$.
\end{lem}

\proof Since $PGL(F,n)$ is the quotient group of $GL(F,n)$ by its
center and the alternating group has trivial  center, Lemma
\ref{dim2} follows from Lemma \ref{dim}. \qed

\subsection{Analysis of Case $I_{32}$}\label{I23}
Since $N$ is a normal subgroup and $N\cap S_d=A_d$, the subgroup
$\langle N,\sigma\rangle$ generated in $G$ by the elements of $N$
and a transposition $\sigma \in S_d$ is isomorphic to a semi-direct
product $N\rtimes \langle \sigma\rangle$. The group $S_d$ is
contained in $\langle N,\sigma\rangle$, since $A_d\subset N$. Thus,
the minimality property of $G$ implies $G=N\rtimes \langle
\sigma\rangle$.

Let us remind that, by assumption, $A_d$ is not a normal subgroup
of $G$.

Suppose first that the group $N$ is not simple.
Pick a minimal non-trivial normal subgroup $N_1$ of $N$. Then, either
$N_1\cap A_d=\{ 1\}$, or $N_1\cap A_d=A_d$, since $A_d$ is simple.

If $N_1\cap A_d=\{ 1\}$, then the group $N_2=\sigma^{-1}N_1\sigma$
is a normal subgroup of $N$ and $N_2\cap A_d=\{ 1\}$. If $N_1=
N_2$, then $N_1$ is a normal subgroup of $G$ and this case is
already considered (Case $I_2$). If $N_1\neq N_2$, then
$[N_1,N_2]\subset N_1\cap N_2=\{ 1\}$ and the group $N_1N_2\simeq
N_1\times N_2$ is a normal subgroup of $G$. Again, the case when
$N_1N_2\cap A_d=\{ 1\}$ is contained in Case $I_2$. Therefore, we
can assume that $N=N_1N_2$. Since, $N_i\cap A_d=\{ 1\}$ for
$i=1,2$, the projections of $A_d$ to the factors should be
imbeddings. Therefore, $|N_i|\geq |A_d|=\frac{d!}{2}$. Hence,
$$k=(G:S_d)=(N:A_d)\geq\frac{d!}{2}$$
which is impossible since, by assumption, $d\geq 1000$ and $k$ is not greater
than either $63d$ or $(4\cdot 42)^2d^2$.

If $N_1\cap A_d=A_d$, then $N_2\cap A_d=A_d$, where
$N_2=\sigma^{-1}N_1\sigma$. If $N_1= N_2$, then $N_1$ is a normal
subgroup of $G$ and this case is contained in Case $I_2$. If
$N_1\neq N_2$, then $N_1\cap N_2$ is a normal subgroup of $N$ and
$A_d\subset N_1\cap N_2\subset N_1$. Therefore, contrary to our
initial assumptions, $N_1$ is not a minimal non-trivial normal
subgroup of $N$.

Thus, it us remains to treat the case when $N$ is a simple group
and $G=N\rtimes \langle\,\sigma\rangle$. Obviously, $N$ can not be
a cyclic group.

If $N$ is isomorphic to some alternating group $A_{d_1}$, then
$d_1-d= n_1\geq 1$ and
$$k=(G:S_d)=(A_{d_1}:A_d)=(d+1)\dots (d+
n_1).$$ By the hypotheses, $d\geq \max(2c, 1000)$ and $k\leq cd^n$,
where either {\rm (i)} $c =63$ and $n=1$, or {\rm (ii)} $c =(4\cdot
42)^2$ and $n=2$. Therefore $n_1\leq 1$ under assumption {\rm (i)}
and $n_1\leq 2$ under assumption {\rm (ii)}.

If $n_1=1$, then $G=S_{d+1}$ (and, moreover, the imbedding of $S_d$
into $G=S_{d+1}$ is the standard one).

Let us show (before ending with an analysis of other simple groups)
that it is impossible to have $n_1=2$ under assumption {\rm (ii)}.

\begin{lem} \label{inAd}
An imbedding $\alpha :A_d \to A_{d+2}$ is conjugate to the standard
one if $d\geq 9$.
\end{lem}

\proof Consider the standard actions of $A_{d}\subset S_{d}$ and
$A_{d+2}\subset S_{d+2}$ on the sets $I_d=\{ 1,2,\dots,d\}$ and
$I_{d+2}=\{ 1,2,\dots,d+2\}$, respectively. If $\tau\in A_d$ is a
cyclic permutation of length $3$, then its image $\alpha(\tau)$ is a
product $\tau_1\dots\tau_s$ of pairwise disjoint cyclic
permutations, and for each $i=1,\dots, s$ it holds $\tau_i^3=1$. To
prove Lemma \ref{inAd}, it suffices to show that $s=1$ for any
$3$-cycle $\tau\in A_d$. Without loss of generality we may assume
that $\tau=(d-2, d-1,d)$.

Under the action of $\alpha(\tau)$, the set $I_{d+2}$ splits into
a disjoint union of $s$ orbits $O_{3,i}$, $i=1,\dots, s$, of
cardinality $3$ and $d+2-3s$ orbits $O_{1,i}$, $i=1,\dots,
d+2-3s$, of cardinality $1$. Consider the subgroup $A_{d-3}$ of
$A_d$ which leaves fixed the elements $d-2,d-1,d\in I_d$. Each
element of $A_{d-3}$ commutes with $\tau$. Hence, the group
$\alpha (A_{d-3})$ acts on the set of orbits $O_{3,i}$ and  on the
set of orbits $O_{1,i}$. This action defines a homomorphism $\beta
: A_{d-3}\to S_s\times S_{d+2-3s}$. But,  $s<d-3$ if $d>9$, and if
$s>1$, then $d+2-3s<d-3$. Therefore, $\beta$ is the trivial
homomorphism if $s>1$, since $A_{d-3}$ is a simple group and
$|A_{d-3}|> |S_s|$, $|A_{d-3}|> |S_{d+2-3s}|$ if $d>9$. Hence, the
homomorphism $\alpha$ induces a homomorphism of $A_{d-3}$ to the
direct product of $s$ copies of $S_3$, which again should be
trivial. Finally, if $s>1$, then $\alpha$ would not be an
imbedding. \qed

Due to Lemma \ref{inAd}, we can assume now that $G\simeq
A_{d+2}\rtimes \langle\,\sigma\rangle$ and that the imbedding
$A_d\subset A_{d+2}$ is standard. Let us show that in this case
$G\simeq  S_{d+2}$ and the imbedding $S_d\subset G\simeq  S_{d+2}$
is standard too.

Indeed, let look at the natural homomorphism $i: Inn(G)\to
Aut(A_{d+2})\simeq S_{d+2}$. Obviously, $i(A_{d+2})=A_{d+2}\subset
S_{d+2}$, and to prove that $G\simeq  S_{d+2}$, it suffices to
show that $i(\sigma)$ is not an inner automorphism of $A_{d+2}$
(recall that $\sigma$ is a transposition as an element of $S_d$).
If $i(\sigma)\in Inn(A_{d+2})$, then there is an element $\tau\in
A_{d+2}$ such that $\gamma=\sigma\tau$ commutes with all the
elements of $A_{d+2}$. In particular, it commutes with $\tau$ and
therefore it commutes with $\sigma$. Since $\sigma\not\in
A_{d+2}$, we have $\gamma=\sigma\tau\neq 1$ and the group $\langle
S_d,\gamma\rangle$ generated in $G$ by $\gamma$ and the elements
of $S_d$ is isomorphic to $S_d\times \langle\gamma\rangle$. But
existence of such a subgroup in $G\simeq A_{d+2}\rtimes
\langle\,\sigma\rangle$ contradicts the minimality property of
$G$.

To show that the imbedding $j:S_d\subset G\simeq  S_{d+2}$ is
standard, note that $j(\sigma)$ is a product $\sigma_1\dots\sigma_s$
of odd number of pairwise disjoint transpositions $\sigma_i\in
S_{d+2}$. We must show that $s=1$. Assume that $s\geq 3$. As in the
proof of Lemma  \ref{inAd}, consider the standard action of $S_{d}$
and $S_{d+2}$ on the sets $I_d=\{ 1,2,\dots,d\}$ and $I_{d+2}=\{
1,2,\dots,d+2\}$, respectively. Let $\sigma\in S_d$ be the
transposition $(d-1,d)$. Under the action of $j(\sigma)$, the set
$I_{d+2}$ splits into a disjoint union of $s$ orbits $O_{2,l}$,
$l=1,\dots, s$, of cardinality $2$ and $d+2-2s$ orbits $O_{1,l}$,
$l=1,\dots, d+2-2s$, of cardinality $1$. Consider the subgroup
$S_{d-2}$ of $S_d$ which leaves fixed the elements $d-1,d\in I_d$.
Each element of $S_{d-2}$ commutes with $\sigma $. Hence, the group
$j(S_{d-2})$ acts on the set of orbits $O_{2,j}$ and  on the set of
orbits $O_{1,j}$. Thus action defines a homomorphism $\beta :
S_{d-2}\to S_s\times S_{d+2-2s}$. But,  $s<d-2$ (recall that $d\geq
1000$), and if $s\geq 3$, then $d+2-2s<d-3$. Therefore the
composition of $\beta$ with the projection to each factor has a
non-trivial kernel if $s\geq 3$. This kernel is either $A_{d-2}$ or
the whole $S_{d-2}$. Therefore the image of each element of
$S_{d-2}$ under the imbedding $j$ has the order not greater than
$4$, which is impossible if $d-2\geq 5$.

The following Lemma forbids an appearance of other simple groups
$N$ in $G=N\rtimes \langle\,\sigma\rangle$ and thus completes the
investigation of Case $I_{23}$.

\begin{lem} \label{index} Assume that a simple group $G$ distinct
from an alternating group contains a subgroup $H_d$ isomorphic
either to the symmetric group $S_d$ or the alternating group
$A_d$, $d\geq 1000$. Then $(G:H_d)>168^2d^2$.
\end{lem}

\proof To prove Lemma, we use the classification of finite simple
groups (see \cite{G}).

The group $G$ is non-abelian, since $H_d$ is a non-abelian group.

The group $G$ can not be a sporadic simple group, since the order
of each sporadic simple group is not divisible by $\frac{d!}{2}$
if $d\geq 33$ (for the sporadic simple groups either the
multiplicity of the prime number $11$ in its order is not greater
than $2$, or the order is not divisible by $13$), while $|G|$ is
divisible by $|H_d|$, which is divisible by $\frac{d!}{2}$.

Let $G$ be a group of Lie type. Then $G$ is a subgroup of either
$GL(F,n)$ or $PGL(F,n)$, where $F$ is a finite field. Denote by
$q$ the number of elements of the field $F$. Since $H_d\subset G$,
then, by Lemmas \ref{dim} and \ref{dim2}, we have the inequality
$n\geq [\frac{d}{4}]$.

If $G$ is one of the following groups: $A_n(q)$,  $B_n(q)$,
$C_n(q)$, $D_n(q)$, $^2A_n(q^2)$,  $^2D_n(q^2)$, then
$$\displaystyle |G|\geq q^{r^2/2}
,$$ where $r=[\frac{n}{2}]$. Since $2^{\frac12
[\frac{d}{8}]^2-16}>(d+2)\log_2 d$ for $d\geq 1000$, we have that
$(G:H_d)>168^2d^2$.

To complete the proof of Lemma, note that all the other simple
groups of Lie type have a non-trivial irreducible linear
representation of dimension less than $250$ and therefore, by Lemma
\ref{dim2}, they can not have a subgroup isomorphic to $A_d$ if
$d\geq 1000$. \qed

\subsection{Analysis of Case $II$}\label{I2}

It follows from Lemma \ref{index} that it remains to consider only
the case $S_d\subset G=A_{d_1}$.

The imbedding $S_d\subset A_{d_1}$ induces an imbedding $S_d\subset
S_{d_1}$. Since any imbedding $S_d\subset S_{d+1}$ is standard, we
have $d_1-d= n_1\geq 2$. By the hypotheses, $d\geq \max(2c, 1000)$
and $(A_{d_1}:S_d)=\frac{1}{2}(d+1)\dots(d+n_1)\leq cd^n$, where
either {\rm (i)} $c =63$ and $n=1$, or {\rm (ii)} $c =(4\cdot 42)^2$
and $n=2$. Therefore $n_1\leq 2$.

Let us show that if $n_1=2$ then the imbedding $S_d\subset A_{d+2}$
is a standard one. Indeed, by Lemma \ref{inAd}, the  imbedding
$S_d\subset A_{d+2}$ induces a standard imbedding $A_d\subset
A_{d+2}$. Moreover, the image in $A_{d+2}\subset S_{d+2}$ of a
transposition $\sigma\in S_d$ is a product of an even number $s$ of
mutually commuting transpositions $\sigma_i$ of $S_{d+2}$. To show
that the imbedding $S_d\subset A_{d+2}$ is a standard one, it
suffices to prove that $s=2$. We omit this proof, since it almost
word by word coincides with the proof of Lemma \ref{inAd}. \qed

\section{Proof of Theorem \ref{main1}.}
\subsection{Minimal expansions of the Galois groups of generic
coverings}\label{3cases} Denote by $\overline g=g-1$ the arithmetic
genus of $X$, $\overline g\ge 1$, and by $B$ the branch locus of
$f:X\to \mathbb P^1$. Since $f$ is $m$-canonical,
$$d=\deg f=2m\overline g.$$
By Hurwitz formula applied to $f$,
$$
\vert B\vert =2d+2\overline g=2(2m+1)\overline g.$$ The branch locus
of $\widetilde f$ (the Galois expansion of $f$) coincides with $B$,
and the ramification indices of the ramification points of
$\widetilde f$ are all equal to $2$. Therefore, by Hurwitz formula
applied to $\widetilde f$,
\begin{equation} \label{widetildeg}
2\widetilde g=-2d!+\frac{1}{2}d!|B|=d!(d+\overline g-2),
\end{equation}
where $\widetilde g=g(\widetilde X)-1$ is the arithmetic genus of
$\widetilde X$.

Assume that $Aut(\widetilde X)\neq Gal(\widetilde X/\mathbb P^1)$
and choose a subgroup $G$ of $Aut(\widetilde X)$ such that
$S_d\subset G$ is a minimal expansion of $S_d$.

Denote by $k=(G:S_d)$ the index of $S_d$ in $G$. The Hurwitz bound
on the order of the automorphism groups of algebraic curves (see,
for example, \cite{FK})  implies that $|G|\leq 84\widetilde g$.
Therefore,
\begin{equation} \label{ineq} k\leq 42(d+\overline g-2).
\end{equation}
In particular, we have
\begin{equation}\label{63} k<63d.
\end{equation}
By Proposition \ref{prop1}, it follows that $G$ is one of the
following groups:
\begin{itemize}
\item[($1$)] $G=S_d\times \mathbb Z/p\mathbb Z$, $p\ge 2$,
$p$ is a prime number;
\item[($2$)] $G=A_d\rtimes D_r$, where
$r\geq 3$, $r$ is odd, $D_r$ is the dihedral group
given by representation
$$
D_r=\langle \sigma, \tau \, |\, \sigma^2=\tau^2=(\sigma\tau)^r=1\rangle,
$$
and the action (by conjugation) of $\sigma$ and $\tau$ on $A_d$
coincides with the action of the transposition $(1,2)\in S_d$ on
$A_d\subset S_d$; \item[($3$)] $G=S_{d+1}$ is the symmetric group.
\end{itemize}

\subsection{Elimination of the remaining three
cases}\label{casebycase} Consider Case $(1)$. Denote by $\gamma$ a
generator of $\mathbb Z/p\mathbb Z$.

Since the action of $
\gamma$ on $\widetilde X$ commutes with the action of any element of
$S_d$, the action of the group $\langle
\gamma\rangle$ on $\widetilde X$ descends to both $X$ and $\mathbb
P^1$.
Denote by $\widetilde X_1=\widetilde X/\langle
\gamma\rangle$, $X_1=X/\langle
\gamma\rangle$ the corresponding factor-spaces and by $\widetilde
r:\widetilde X\to \widetilde X_1$, $r:X\to X_1$, $h_1:\widetilde
X_1\to X_1$, and $r_P:\mathbb P^1\to\mathbb P^1/\langle
\gamma\rangle\simeq \mathbb P^1$ the corresponding morphisms. We
have the following commutative diagram:

$$
\xymatrix{ \widetilde X \ar[d]_{\widetilde r} \ar[r]^{h } & X
\ar[d]^{r} \ar[r]^{f } & \mathbb P^1 \ar[d]^{
r_P} \\
\widetilde X_1 \ar[r]_{h_1 } & X_1
\ar[r]_{f_1 } & \mathbb P^1. }
$$

The cyclic covering $ r_P:\mathbb P^1\to \mathbb P^1$ is of degree
$p\ge 2$ and it is ramified at two points, say $x_1,\, x_2\in
\mathbb P^1$. Therefore, the cyclic covering $r$ is ramified at
least at $2(d-1)$ points lying in $f^{-1}(x_1)\cup f^{-1}(x_2)$ and
their ramification index is equal to $p$. By Hurwitz formula,
$$ 2\overline g\geq 2p(g(X_1)-1)+2(d-1)(p-1)
$$
which implies
$$2\overline g\geq -2p+2(2m\overline g-1)(p-1).$$
Finally, thus we  get the inequality
$$p\leq \frac{(2m+1)\overline g+1}{2m\overline g-
2}=1+ \frac{\overline g+3} {2m\overline g-2} < 2$$ which shows that
Case $(1)$ is impossible.

Consider Case $(2)$. For a suitable pair of generators $\sigma$,
$\tau$ of $D_r$, we have $S_d=A_d\rtimes \langle
\sigma\rangle\subset G$, while the group $S^{\prime}_d=A_d\rtimes
\langle \tau\rangle$ is conjugated to $S_d$ and does not coincide
with $S_d$ (but, it is isomorphic to $S_d$). Besides, $A_d\subset
S^{\prime}_d\cap S_d$. Denote by $X_1=\widetilde
X/S^{\prime}_{d-1}$, $\mathbb P^1=\widetilde X/S^{\prime}_{d}$, and
$X_0=\widetilde X/A_{d}$ the corresponding quotient spaces. They can
be arranged in the following commutative diagram in which the
morphisms $f_{0i}$, $i=1,2$, are of degree two and, since $f$ is a
generic covering, $f_{01}$ is branched over all the points belonging
to $B$.

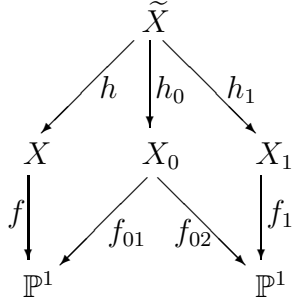
\begin{picture}(300,125)
\put(170,110){$\widetilde X$} \put(125,60){$X$} \put(170,60){$X_0$}
\put(213,60){$X_1$} \put(125,10){$\mathbb P^1$}
\put(213,10){$\mathbb P^1$}
\put(167,107){\vector(-1,-1){35}}\put(178,107){\vector(1,-1){35}}
\put(173,107){\vector(0,-1){33}}
\put(176,55){\vector(1,-1){33}}\put(172,55){\vector(-1,-1){33}}
\put(127,56){\vector(0,-1){32}} \put(215,56){\vector(0,-1){32}}
\put(154,85){$h$}\put(175,85){$h_0$}\put(202,85){$h_1$}
\put(119,38){$f$}\put(157,31){$f_{01}$} \put(182,31){$f_{02}$}
\put(217,38){$f_1$}
\end{picture}
\newline The degree 2 morphisms  $f_{0i}$, $i=1,2$, define an imbedding
$i: X_0\to \mathbb P^1\times \mathbb P^1$ with $i(X_0)$ being a
curve of bi-degree $(2,2)$ in $\mathbb P^1\times \mathbb P^1$.
Therefore, $i(X_0)$ is an elliptic curve and the projection of
$i(X_0)$ onto each factor is branched at four points. On the other
hand, $f_{01}$ is branched at every point of $B$ and $\vert
B\vert=2d+2\overline g>4$. Therefore,  Case $(2)$ is impossible.

Consider Case $(3)$. Note, first of all, that the imbedding of $S_d$
into $G=S_{d+1}$ is the standard one.

Consider the quotient space $\widetilde X/G$ and the quotient map
$\overline f: \widetilde X\to \widetilde X/G$. The latter factors
through $\widetilde f$, so that $\widetilde X/G\simeq \mathbb P^1$
and $\overline f$ is the composition of the following morphisms

$$
\xymatrix{
\widetilde X  \ar[r]^{h } & X
\ar[r]^{f } & \mathbb P^1 \ar[r]^{r } & \mathbb P^1,  }
$$
where $r$ is a morphism of degree $d+1$. Since $S_{d}$ and $S_{d+1}$
have no common normal subgroups, $\overline f$ is the Galois
expansion of $r$.

Denote by $B_1\subset \mathbb P^1$ the branch locus of $r$ and
compare the cardinality of $B$ with the cardinality of $r(B)\subset
B_1$.

The symmetric group $S_{d+1}$ acts as the permutation group on the
set $I=\{ 1,\dots, d+1\}\subset \mathbb N$. Denote by $H_i=\{ \gamma
\in S_{d+1}\, |\, \gamma(i)=i\}$, so that our $S_d=H_{d+1}$. All
groups $H_i$ are conjugated to each other. Therefore for each $i$
the covering $\widetilde f_i: \widetilde X\to \widetilde X/H_i\simeq
\mathbb P^1$ is the Galois expansion of a generic covering.

Let $a\in \widetilde X$ be a ramification point of $\widetilde
f_{d+1}=\widetilde f$. The stabilizer group $St_a(\overline f)=\{
g\in G\, |\, g(a)=a\}$ is a cyclic group; its order is equal to the
ramification index of $\overline f$ at $a$. Let $\tau$ be a
generator of $St_a(\overline f)$. The intersection $St_a(\overline
f)\cap S_d=St_a(\widetilde f_{d+1})$ is a group of order two
generated by a transposition $\sigma\in H_{d+1}$, since $\widetilde
f_{d+1}$ is the Galois expansion of a generic covering. Therefore,
$\sigma =\tau^k$, where $\tau^{2k}=1$,

Let us show, first, that $k$ is odd. Indeed, let us write $\tau$
as the product of cyclic permutations: $\tau=(i_{1,1},\dots,
i_{1,k_1})\dots (i_{s,1},\dots, i_{s,k_s})$. We can assume that,
up to renumbering, $\sigma=(i_{1,1},\dots, i_{1,k_1})^k$ and
$(i_{j,1},\dots, i_{j,k_j})^k=1$ for $j=2,\dots, s$. Now, it easy
to see that $k_1=2$, $k$ is odd, and all $k_j$ are divisors of $k$
for $j=2,\dots, s$.

Let us show that $k=1$, so that $\tau =\sigma$. We have for each $i$
the intersection $St_a(\overline f)\cap H_i=St_a(\widetilde f_{i})$
is a group of order at most two and if its order equal to two, then
it is also generated by a transposition $\sigma_i\in H_{i}$, since
$\widetilde f_{i}$ is the Galois expansion conjugated to $\widetilde
f_{d+1}$. On the other hand, the element
$$\sigma\tau=(i_{1,2},\dots, i_{1,k_2})\dots (i_{s,1},\dots, i_{s,k_s})\in
St_a(\overline f)\cap H_{i_{1,1}}=St_a(\widetilde f_{i_{1,1}})$$
is of odd order. Therefore  $\sigma=\tau$.

Now, consider the fibre $\overline f^{-1}(\overline f(a))$
containing the point $a$. The fibre $\overline f^{-1}(\overline
f(a))$ can be identified with the set of right  cosets
$\{St_a(\overline f)\gamma \}$ in $S_{d+1}$. The stabilizer group
$St_{(a)\gamma}(\overline f)$ of the point $(a)\gamma$ is generated
by the transposition $\gamma^{-1}\sigma\gamma$.

The fibre $\overline f^{-1}(\overline f(a))$ splits into the
disjoint union of orbits under the action of $H_{d+1}$ each of which
is a fibre of $\widetilde f_{d+1}$. Without loss of generality, we
can assume that the group $St_a(\overline f)$ is generated by
$\sigma=(1,2)$. Then it is easy to see that each of these orbits can
be identified with one of $F_i=\{St_a(\overline f)\gamma \, |\,
\gamma \in \sigma_iH_{d+1}\}$, where $\sigma_i=(i,d+1)$ if $2\leq
i<d+1$, and $\sigma=(1,2)$, if $i=d+1$. (The transpositions
$\sigma_1=(1,d+1)$ and $\sigma_2=(2,d+1)$ give the same orbit under
the action of $H_{d+1}$, since $(1,2)(1,d+1)(1,2)= Id\cdot
(2,d+1)$.)

Now, the points
$a_i=(a)
\sigma_i$ have the same stabilizer group
$$St_{a_i}(\overline f)=\langle (1,2)\rangle\subset H_{d+1}$$ if $i>2$.
Therefore for $i\geq 3$ the points belonging to $F_i$ are the
ramification points of $\widetilde f_{d+1}$ and hence  $\widetilde
f_{d+1}(F_i)\in B$. It is easy to see that the points belonging to
$F_2$ are not the ramification points of $\widetilde f_{d+1}$.
Therefore the point $\widetilde f_{d+1}(F_2)$ is a ramification
point of $r$.

As a consequence, we obtain that if $\widetilde b\in r(B)$, then the
fibre $r^{-1}(\widetilde b)$ consists of $d-1$ points belonging to
$B$ and one point (a ramification point of $r$) which does not
belong to $B$. Hence, $|B|=2d+2\overline g=2(m+1)\overline g$ is
divisible by $d-1=2m\overline g-1$. Then, $2\overline g+1$ should
also be divisible by $2m\overline g-1$. But, it is possible only if
$\overline g=1$ and $m=1$ or $2$. \qed

\section{Proof of Theorem \ref{main2}.}

\subsection{Local behavior of generic coverings and their Galois
expansions.} \label{sec31} In this subsection we specialize to
surfaces the definitions related to generic actions of the symmetric
group, compare our definitions with the traditional definition of
generic coverings, deduce the local behavior of generic coverings
from the local behavior of these actions, and fix the corresponding
notation and notions.

Recall that the Galois expansion $\widetilde f:\widetilde X\to
\mathbb P^2$ of a generic covering $f: X\to \mathbb P^2$ of degree
$d$ is factorized through $f$, $\widetilde f=f\circ h$, by means of
a Galois covering $h:\widetilde X\to X$ with the Galois group
$Gal(\widetilde X/X)=S_{d-1}\subset S_d=Gal(\widetilde X/\mathbb
P^2)$. The branch locus $B\subset \mathbb P^2$ of $f$ coincides with
that of $\widetilde f$. We have called $f$ generic if the action of
$S_d$ on $\widetilde X$ is generic. The latter means that for any
point $a\in \widetilde X$ its stabilizer $St_a(S_d)$ is a standard
imbedded in $S_d$ product of symmetric groups (depending on $a$) and
the actions induced by $St_a(S_d)$ on the tangent spaces
$T_a\widetilde X$ are standard representations of rank $\leq 2$ (see
Introduction).

On the other hand, in dimension two more traditionally one
understands under a generic covering $f:X\to \mathbb P^2$ of
degree $d$ a covering whose local behavior is as follows. The
branch locus $B$ of $f$ is a cuspidal curve. Over a neighborhood
$U$ of a smooth point of $B$ the preimage $f^{-1}(U)$ splits into
a disjoint union of $d-1$ connected components, in one of them the
covering is two-sheeted and isomorphic to the projection to $x,y$
plane of the surface $x=z^2$ at a neighborhood of
the origin,
and in the other components it is a local isomorphism. Over a
neighborhood of a cuspidal point of $B$ the preimage splits into a
disjoint union of $d-2$ neighborhoods, in one of them the covering
is a pleat, that is a three-sheeted covering which is isomorphic to
the projection to $x,y$ plane of the surface $y=z^3+xz$ at a
neighborhood of the origin,
and in the other components it is a local isomorphism. Over a
neighborhood of a node of $B$ the preimage splits into a disjoint
union of $d-2$ neighborhoods, in two of them this covering is
two-sheeted and isomorphic in the union of them to the projection to
$x,y$ plane of two surfaces $x=z^2$ and $y=z^2$ in a neighborhood of
the origin and in the other components it is a local isomorphism.
The nontrivial local Galois groups in the corresponding three cases
are $\mathbb Z/2$, $S_3$, and $\mathbb Z/2\times\mathbb Z/2$. But
all their nontrivial representations in $GL(2,\mathbb C)$ which
produce a non-singular quotient are standard representations of rank
$\leq 2$, and therefore our definition of generic coverings
coincides with the traditional one.

The local behavior of generic coverings is easily understandable from the above
local models. In particular, one observes that: $f^*(B)=2R+C$, where $R$ (the
ramification locus of $f$) is nonsingular, $C$ is reduced and non-singular above the non-singular
points of $B$; $R$ and $C$ intersect each other only
above the nodes and cusps of $B$; they meet at two points above each node and intersect there
transversally; and they meet at one point above each cusp and intersect there with simple
tangency.

Let us observe the same local behavior from the point of view of the action of $S_d$
on $\widetilde X$, which will help us in our further considerations in this section.

At a small neighborhood of any point $a\in\widetilde X$ the action
of $St_a(S_d)$ can be linearized (Cartan's linearization procedure
\cite{Cartan}, which by a suitable change of coordinates identifies
the local action with the action induced on the tangent space
$T_a\widetilde X$, is reproduced  below in the proof of Lemma
\ref{stability}). Let us treat case by case different possibilities
for $St_a(S_d)$ and, by means of a linearization of the action and
in accordance with the definition of standard representations,
analyze the local behavior of $f$, $\widetilde f$, and $h$.

If $St_a(S_d)=S_2$ is generated by a transposition $\sigma\in S_d$,
then in a neighborhood of $a$ the ramification locus $\widetilde R$
of $\widetilde f$ coincides with the set of fixed points of
$\sigma$, which we denote by $\widetilde R_\sigma $. The latter is
smooth everywhere, and, in particular, $\widetilde R$ is smooth at
$a$. The image $h(a)$ of $a$ belongs to the ramification locus $R$
of $f$ (equivalently, $a$ does not belong to the ramification locus
of $h$) if and only if $\sigma \not\in S_{d-1}$. Moreover,
$h(\widetilde R_\sigma)$ coincides with $R$ at a neighborhood of
$h(a)$ if $\sigma \not\in S_{d-1}$ (otherwise, it coincides with $C$
introduced above). In both cases, $\sigma \in S_{d-1}$ and $\sigma
\not\in S_{d-1}$, the curve $h(\widetilde R_\sigma)$ is smooth at
$h(a)$. Furthermore, in both cases, $\widetilde f(a)$ belongs to $B$
and $B$ is non-singular at $\widetilde f(a)$.

If $St_a(S_d)=S_2\times S_2$ is generated by two commuting
transpositions, $\sigma_1\in S_d$ and $\sigma_2\in S_d$, then the
point $a$ belongs to $\widetilde R_{\sigma_1}\cap \widetilde
R_{\sigma_2}$, the curves $\widetilde R_{\sigma_1}$ and $\widetilde
R_{\sigma_2}$ are nonsingular, and they meet transversally at $a$.
Furthermore, $h(\widetilde R_{\sigma_1})$ and $h(\widetilde
R_{\sigma_2})$ are nonsingular and meet transversally. If one of the
transpositions, say $\sigma_1$, does not belong to $S_{d-1}$ the
curve $h(\widetilde R_{\sigma_1})$ is contained in $R$  and,
moreover, coincides with $R$ in a neighborhood of $h(a)$. If
$\sigma_1\in S_{d-1}$ the curve $h(\widetilde R_{\sigma_1})$ is not
contained in $R$ (but contained in $C$). If both $\sigma_1$ and
$\sigma_2$ belong to $S_{d-1}$, $h(a)$ is not a ramification point
of $f$ (and then it is a node of $C$ with $C=h(\widetilde
R_{\sigma_1})\cup h(\widetilde R_{\sigma_2})$ in a neighborhood of
$h(a)$). In any case, $\widetilde f(a)$ is a node of $B$.

If $St_a(S_d)=S_3$ is generated by two non-commuting transpositions,
$\sigma_1\in S_d$ and $\sigma_2\in S_d$, then the point $a$ belongs
to $\widetilde R_{\sigma_1}\cap \widetilde R_{\sigma_2}\cap
\widetilde R_{\sigma_3}$, where $\sigma_3=\sigma_1\sigma_2\sigma_1$,
the curves $\widetilde R_{\sigma_1}$,  $\widetilde R_{\sigma_2}$,
and $\widetilde R_{\sigma_3}$ are nonsingular and meet pairwise
transversally at $a$. If all the three transpositions belong to
$S_{d-1}$, the point $h(a)$ is not a ramification point of $f$.
Otherwise, one and only one of the transpositions, say $\sigma_3$,
belongs to $S_{d-1}$, and then: $h(\widetilde
R_{\sigma_1})=h(\widetilde R_{\sigma_2})$ and $h(\widetilde
R_{\sigma_3})$ are nonsingular, they are tangent to each other, and
$h(\widetilde R_{\sigma_1})=h(\widetilde R_{\sigma_2})$ coincides
with $R$ (while $h(\widetilde R_{\sigma_3})$ coincides with $C$) in
a neighborhood of $h(a)$. In any case, $\widetilde f(a)$ is a cusp
of $B$.

\subsection{Invariants of $m$-canonical generic coverings.}\label{sec32}

Assume that $f:X\to \mathbb P^2$ is a generic $m$-canonical covering
branched along a cuspidal curve $B\in \mathbb P^2$. Then $X$ is a
minimal surface of general type, it does not contain any
$(-2)$-curve, and the degree of $f$ is equal to
$$d=\deg f=m^2K_{X}^2.$$
According to the formula for the canonical divisor of a finite
covering, $K_X=f^*K_{\mathbb P^2}+[R]$. Hence, the divisor $R$ is
numerically equivalent to $(3m+1)K_X$. Since, in addition, the curve
$R$ is non-singular and $X$ has no $(-2)$-curves (if $f$ is a
$m$-canonical generic covering, then $K_X$ is ample), $R$ is
irreducible. Therefore, $B$ as a curve birational to $R$ is
irreducible as well. Thus, we can apply the results from \cite{Ku}.
In particular, we get the following formulas for the degree $\deg B$
and the number $c$ of cusps of $B$ (see the proof of Theorem 2 in
\cite{Ku}):
\begin{equation} \label{degB}
\deg B =m(3m+1)K^2_X
\end{equation}
and
\begin{equation} \label{c}
c =(12m^2+9m+3)K^2_X -12p_a, \end{equation} where $p_a=p_g-q+1$ is
the arithmetic genus of $X$. Note that if $\deg f\geq 3$ for a
generic covering $f$, then its branch curve $B$ should have cuspidal
singular points, that is, $c>0$ (indeed, the image in $S_d$ of the
monodromy of $f$ is a transitive subgroup of $S_d$ and, for generic
coverings, this image is generated by transpositions, hence
coincides with $S_d$; therefore $\pi_i(\mathbb P^2\setminus B)$ is
non-abelian if $d\ge 3$, while by Zariski's theorem $\pi_i(\mathbb
P^2\setminus B)$ is an abelian group if $B$ is a nodal curve).

Finally, applying the projection formula for the canonical divisor
to $\widetilde f$ we obtain $K_{\widetilde X}=\widetilde
f^*(K_{\mathbb P^2})+[\widetilde R]=\widetilde f^*(K_{\mathbb
P^2}+\frac12[B])$, and therefore
\begin{equation} \label{widetilKX}
K^2_{\widetilde X}=\frac{1}{4}(\deg B -6)^2d!
=d!(\frac{m(3m+1)}2K^2_X-3)^2.
\end{equation}

\subsection{Minimal expansions of the Galois groups of generic
coverings}\label{4cases} Assume that $Aut(\widetilde X)\neq
Gal(\widetilde X/\mathbb P^2)$ and choose a subgroup $G$ of
$Aut(\widetilde X)$ such that $S_d\subset G$ is a minimal expansion
of $S_d$. Denote by $k=(G:S_d)$ the index of $S_d$ in $G$.

The Xiao bound (see \cite{Xi}) on the order of the automorphism
groups of surfaces of general type states that $|G|\leq
42^2K^2_{\widetilde X}$. It implies
\begin{equation} \label{ineq2} k\leq 42^2(\frac{m(3m+1)}2K^2_X-3)^2.
\end{equation}
Finally we get
\begin{equation}\label{63??}
k<(2\cdot 42)^2d^2.
\end{equation}
By Proposition  \ref{prop1}, it follows that the group $G$ can be
only one of the following groups:
\begin{itemize}
\item[($1$)] $G=S_d\times \mathbb Z/p\mathbb Z$, $p\ge 2$,
$p$ is a prime number;
\item[($2$)] $G=A_d\rtimes D_r$, where
$r\geq 3$, $r$ is odd, $D_r$ is the dihedral group
given by representation
$$
D_r=\langle \sigma, \tau \, |\, \sigma^2=\tau^2=(\sigma\tau)^r=1\rangle,
$$
and the action (by conjugation) of $\sigma$ and $\tau$ on $A_d$
coincides with the action of the transposition $(1,2)\in S_d$ on
$A_d\subset S_d$;
\item[($3$)] $G=S_{d+1}$ is the symmetric group;
\item[($4$)] $G=A_{d+2}$ is the alternating group,
and the imbedding of $S_d$ into $G=A_{d+2}$ is a standard one.
\end{itemize}

\subsection{Case $(1)$.} Denote by $g$ a generator of $\mathbb Z/p\mathbb Z$. As
in the proof of Theorem \ref{main1}, since the action of $g$ on
$\widetilde X$ commutes with the action of any element of $S_d$, the
action of the group $\langle g\rangle=\mathbb Z/p\mathbb Z$ on
$\widetilde X$ descends to both $X$ and $\mathbb P^2$. Denote by
$\widetilde X_1=\widetilde X/\langle g\rangle$, $X_1=X/\langle
g\rangle$ the corresponding factor-spaces and by $\widetilde
r:\widetilde X\to \widetilde X_1$, $r:X\to X_1$, $h_1:\widetilde
X_1\to X_1$, and $ r_P:\mathbb P^2\to\mathbb P^2/\langle g\rangle=Y$
the corresponding morphisms. We have the following commutative
diagram:

$$
\xymatrix{ \widetilde X \ar[d]_{\widetilde r} \ar[r]^{h } & X
\ar[d]^{r} \ar[r]^{f } & \mathbb P^2 \ar[d]^{
r_P} \\
\widetilde X_1 \ar[r]_{h_1 } & X_1
\ar[r]_{f_1 } & Y. }
$$

The automorphism $g$ on $\mathbb P^2$ is defined by a linear map
$\mathbb C^3\to\mathbb C^3$ of period $p$. Therefore, it has  either
three isolated fixed points, say $x_1,\, x_2, x_3\in \mathbb P^2$,
or an isolated fixed point, say $x\in \mathbb P^2$ and a fixed line
$E\subset \mathbb P^2$. Respectively, the cyclic covering
$r_P:\mathbb P^2\to \mathbb P^2$ is of degree $p$ and it is ramified
either at three points $x_1,\, x_2, x_3$ or at the point $x$ and
along the line $E$. Consider this two cases separately.

If $r_P$ is ramified at three isolated points, then $r$ is ramified
at at least $3(d-2)$ points lying in $F=f^{-1}(x_1)\cup
f^{-1}(x_2)\cup f^{-1}(x_3)$ and the ramification index of each of
them is equal to $p$. The points of $F$ are the only fixed points of
the automorphism $g$ acting on $X$. Therefore, by Lefcshetz fixed
point theorem we have
\begin{equation} \label{Lef}
\displaystyle |F|=\sum _{i=0}^4(-1)^itr_i,
\end{equation}
where $tr_i$ is the trace of the linear transformation $g^*$ acting
on $H^i(X,\mathbb R)$. It follows from (\ref{Lef}) that
\begin{equation} \label{Lef1}
\displaystyle 3(d-2)\le \vert F\vert\leq  \sum _{i=0}^4|tr_i|\leq
\sum _{i=0}^4b_i(X)=e(X)+4b_1(X)
\end{equation}
(where $e$ stands for the topological Euler characteristic). On the
other hand, Noether's formula $1-q+p_g=\frac{K^2_X+e(X)}{12}$ where
$2q=b_1$ gives us
\begin{equation} \label{Lef2}
\displaystyle e(X)+ 4b_1(X)= 12-K^2_X+12p_g-4q \leq  12-K^2_X+12p_g
\end{equation}
Therefore, combining (\ref{Lef1}) and  (\ref{Lef2}) with Noether's
inequality $2p_g\le K^2_X+4$ we get
\begin{equation} \label{Lef3}
\displaystyle 3(m^2K^2_X-2)\leq 12-K^2_X+12p_g \leq 5K^2_X+36.
\end{equation}
Hence,
$$ (3m^2-5)K_X^2\leq 42 $$
which contradicts to $m^2K_X^2\geq  2\cdot 84^2$ if $m\geq 2$, since
$K_X^2\geq 1$ for any minimal surface of general type.

Now, let us assume that $r_P$ is ramified at a point $x\in\mathbb
P^2$ and along a line $E\subset\mathbb P^2$. In this case each line
$L\subset \mathbb P^2$ passing through $x$ is invariant under the
action of $g$ on $\mathbb P^2$. Therefore, each curve
$C=f^{-1}(L)\subset X$ is invariant under the action of $g$ on $X$.
Pick a generic line $L$ passing through $x$. By Hurwitz formula,
\begin{equation} \label{gen1}
2(g(C)-1)=-2d +\deg B -2m^2K_X^2+m(3m+1)K^2_X
=m(m+1)K^2_X,
\end{equation}
where $g(C)$ is the geometric genus of $C$.

Consider the restriction $r_{\mid C} :C\to C/\langle
g\rangle=Z\subset X_1$ of $r$ to $C$. The cyclic covering $r_{\mid
C}$ has degree $p$ and it is branched at at least $2d-3=2m^2K_X^2-3$ points
belonging to $f^{-1}(x)\cup f^{-1}(L\cap E)$. Therefore,
\begin{equation} \label{gen2} 2(g(C)-
1)\geq 2p(g(Z)-1)+ (2m^2K_X^2-3)(p-1).
\end{equation}
It follows from (\ref{gen1}) and (\ref{gen2}) that
$$m(m+1)K^2_X
\geq -2p+ (2m^2K_X^2-3)(p-1)$$ which implies
$$m(3m+1)K^2_X\geq
(2m^2K_X^2-5)p. $$ Finally, since $(m^2-m)K^2\ge 500>10$, we get the
inequality
$$p\leq \frac{m(3m+1)K^2_X} {2m^2K_X^2-5}<2,
$$
which is a contradiction.

\subsection{Case $(2)$. Group theoretic part}\label{case2part1}
Since $r$, $r\ge 3$, is odd, the conjugacy
class of $\sigma$ in $D_r$ consists of $r$ elements
$\sigma_1=\sigma,\sigma_2=\tau,\dots,\sigma_r$. For each $i, 1\le
i\le r,$ the group $S_{d,i}$ generated in $G$ by $\sigma_i$ and the
elements of $A_d$ is isomorphic to $A_d\rtimes \langle
\sigma_i\rangle\simeq S_d$. The groups $S_{d,i}$ are conjugate to
each other in $G$. Besides, $\sigma_i \in S_{d,i}$ acts on
$A_d\subset S_{d,i}$ as the transposition $(1,2)$. The element of
$S_{d,i}$ which is conjugate to $\sigma_i$ and acts on $A_d$ as a
transposition $(i_1,i_2)$ will be denoted by $\sigma_{i,(i_1,i_2)}$.
Given two disjoint subsets $J_1\ne\emptyset, J_2$ of $I=\{
1,\dots,d\}$, denote by $S_{J_1\bigsqcup J_2,i}$ the subgroup of
$G=A_d\rtimes D_r$ generated by the elements $\sigma_{i,(i_1,i_2)
}$, $(i_1,i_2)\in (J_1\times J_1)\cup (J_2\times J_2)$.

Let $St_a\subset Aut \widetilde X$ be the stabilizer of a point
$a\in \widetilde X$. For a subgroup $H$ of $Aut \widetilde X$ put
$St_a(H)=H\cap St_a$. For each point $a\in \widetilde X$ the action
induced by $St_a(S_d)$ on the tangent space $T_a\widetilde X$ is a
standard representation of rank $\leq 2$, and the group $St_a(S_d)$
is trivial or can be expressed as $S_{J_1\bigsqcup J_2,1}$, where
either $2\le |J_1|\le 3$ and $J_2=\emptyset$, or $|J_1|=|J_2|=2$.
Since the groups $S_{d,i}$ are conjugate to each other, for each $i$
and for each point $a\in \widetilde X$ the group $St_a(S_{d,i})$ has
the same properties. Therefore, the intersection $St_a(S_{d,i})\cap
A_d$ is generated by the cyclic permutation $(i_1,i_2,i_3) \in A_d$
if $J_1=\{ i_1,i_2,i_3\}$ and $J_2=\emptyset$, and it is generated
by the product of two transpositions if $J_1=\{ i_1,i_2\}$ and
$J_2=\{ i_3,i_4\}$. In the remaining cases ($|J_1|=2$ and
$J_2=\emptyset$, or $St_a(S_{d,i})=\{ 1\}$) the group
$St_a(S_{d,i})\cap A_d$ is trivial. This implies that if
$St_a(S_{d,1})=S_{J_1\bigsqcup J_2,1}$, where $|J_1|=3$ and
$J_2=\emptyset$, then $St_a(S_{d,i})=S_{J^{\prime}_1 \bigsqcup
J^{\prime}_2,i}$, where $J^{\prime}_1=J_1$ and
$J^{\prime}_2=\emptyset$, for each $i$. Similarly, if
$St_a(S_{d,1})=S_{J_1\bigsqcup J_2,1}$, where $J_1=\{ i_1,i_2\}$ and
$J_2=\{ i_3,i_4\}$, then $St_a(S_{d,i})=S_{J^{\prime}_1\bigsqcup
J^{\prime}_2,i}$, where $J^{\prime}_1=J_1$ and $J^{\prime}_2=J_2$,
for each $i$. If $St_a(S_{d,1})=S_{J_1\bigsqcup J_2,1}$, where
$J_1=\{ i_1,i_2\}$ and $J_2=\emptyset$, then for each $i$ either
$St_a(S_{d,i})$ is trivial or
$St_a(S_{d,i})=S_{J^{\prime}_1\bigsqcup J^{\prime}_2,i}$, where
$|J^{\prime}_1|=2$ and $J^{\prime}_2=\emptyset$.

Let us examine more in details the case
$St_a(S_{d,i})=S_{J_1\bigsqcup J_2,i}$ for each $i$, where $|J_1|=3$
and $J_2=\emptyset$. Denote by $y$ the cyclic permutation
$(i_1,i_2,i_3)\in A_d$, where $J_1=\{ i_1,i_2,i_3\}$.
The group $St_a(G)$ contains a subgroup $F_{3}$ generated by three
elements $x_1=\sigma_{1,(i_1,i_2)}$ (conjugate to $\sigma$),
$x_2=\sigma_{2,(i_1,i_2)}$ (conjugate to $\tau$), and $y$. It is
easy to see that $F_3$ has the following presentation:
\begin{equation}\label{repr}
\begin{array}{ll}
F_3=\langle x_1,x_2,y\, \mid \, &
x_1^2=x_2^2=(x_1x_2)^r=y^3=[y,x_1x_2]=[y,x_2x_1]=1, \\
&   x_1^{-1}yx_1=y^{-1},\, \,  x_2^{-1}yx_2=y^{-1} \rangle
\end{array}
\end{equation}
(recall that $r\geq 3$). The group $F_3$ is non-abelian. It contains
a maximal normal subgroup $N_3$ generated by $y$ and $z=x_1x_2$.
This subgroup is isomorphic to the direct product $\langle
y\rangle\times \langle z\rangle$ of two cyclic groups of orders $3$
and $r$, respectively. On the other hand, according to well known
properties of finite subgroups of $GL(2, \mathbb C)$ (see, for
example, \cite{Dorn}) the quotient of $F_3$ by its center should be
either a cyclic group, or a dihedral group, or $A_4$, or $A_5$, or
$S_5$. But, $F_3= (\mathbb Z/3\mathbb Z\times\mathbb Z/r\mathbb
Z)\rtimes \mathbb Z/2\mathbb Z$ and thus it has the trivial center,
since $r$ is odd. All together, these arguments imply that $r$
should not be divisible by $3$ (recall that the branch curve does
have at least one cuspidal point, see Subsection \ref{sec32}) and
$F_3$ should be isomorphic to the dihedral group $D_{r^{\prime}}$,
$r^{\prime}=3r$.

In addition, once more due to the known classification of conjugacy
classes of finite subgroups of $GL(2,\mathbb C)$, the action of
$F_3\simeq D_{3r}$ near the point $a$ is isomorphic to the unique
$2$-dimensional linear representation of $D_{3r}$. In particular, at
a neighborhood of $a$ the fixed points of $\sigma_{i,(i_1,i_2)}$,
$i_1,i_2\in J_1$, form a smooth curve, which we denote by
$\widetilde R_{i,(i_1,i_2)}$, and any two of them, $\widetilde
R_{i,(i_1,i_2)}$ and $\widetilde R_{i',(i'_1,i'_2)}$ with
$(i,(i_1,i_2))\ne (i',(i'_1,i'_2))$, are distinct and meet each
other at $a$ transversally.

Next, let us examine the case $St_a(S_{d,i})=S_{J_1\bigsqcup
J_2,i}$, where $J_1=\{ i_1,i_2\}$ and  $J_2=\{ i_3,i_4\}$. The group
$St_a(G)$ contains a subgroup $F_{2,2}$ generated by three elements
$x_1=\sigma_{1,(i_1,i_2)}$ (conjugate to $\sigma$),
$x_2=\sigma_{2,(i_1,i_2)}$ (conjugate to $\tau$), and
$y=(i_1,i_2)(i_3,i_4)\in A_d$. It is easy to see that $F_{2,2}$ has
the following presentation:
\begin{equation}\label{repr2}
F_{2,2}=\langle x_1,x_2,y\, \mid \,
x_1^2=x_2^2=(x_1x_2)^r=y^2=[y,x_1]=[y,x_2]=1 \rangle
\end{equation}
(recall that $r\geq 3$ and $r$ is odd). The group $F_{2,2}$ contains
a maximal normal subgroup $N_{2,2}$ generated by $y$ and $z=x_1x_2$.
This subgroup is isomorphic to the direct product $\langle
y\rangle\times \langle z\rangle$ of two cyclic groups of orders $2$
and $r$, respectively. Therefore, $F_{2,2}$ is isomorphic to the
dihedral group $D_{2r}$. According to the classification of
conjugacy classes of finite subgroups of $GL(2,\mathbb C)$, the
action of $F_{2,2}\simeq D_{2r}$ near the point $a$ is isomorphic to
the unique $2$-dimensional linear representation of $D_{2r}$. In
particular, similar to the previous case, for any
$i_1,i_2,i'_1,i'_2\in J_1\bigsqcup J_2$, the curves $\widetilde
R_{i,(i_1,i_2)}$ and $\widetilde R_{i',(i'_1,i'_2)}$ with
$(i,(i_1,i_2))\ne (i',(i'_1,i'_2))$ are distinct and meet each other
at $a$ transversally.

Finally, consider the case $St_a(S_{d,1})=S_{J_1\bigsqcup J_2,1}$
and $St_a(S_{d,2})=S_{J_1^{\prime}\bigsqcup J_2^{\prime},2}$, where
$J_1=\{ j_1,j_2\}$,  $J_2=\emptyset$, $J^{\prime}_1=\{ j_3,j_4\}$,
$J^{\prime}_2=\emptyset$. Let us show that either $J_1\cap
J^{\prime}_1=\emptyset$ and then $St_a(G)$ contains a subgroup
isomorphic to $D_r$, or $J_1=J^{\prime}_1$.

Indeed, if $|J_1\cap J^{\prime}_1|=1$, then we can assume that
$j_2=j_3$ so that
$\sigma_{2,(j_3,j_4)}=\eta^{-1}\sigma_{2,(j_1,j_2)}\eta$, where
$\eta=(j_4,j_2,j_1)\in A_d$.  We have $\eta^3=1$ and
$$\begin{array}{ll}
(\sigma_{1,(j_1,j_2)}\sigma_{2,(j_3,j_4)})^r= &
(\sigma_{1,(j_1,j_2)}
\eta\sigma_{2,(j_1,j_2)}\eta^{-1})^r=(\sigma_{1,(j_1,j_2)}
\sigma_{2,(j_1,j_2)}\eta)^r= \\
& (\sigma_{1,(j_1,j_2)} \sigma_{2,(j_1,j_2)}
)^r\eta^r=\eta^r=\eta^{\pm 1},
\end{array}
$$
since $r$ is not divisible by $3$. Therefore, $\eta\in St_a(G)$
which contradicts the assumption that $St_a(S_{d,1})=S_{J_1\bigsqcup
J_2,1}$.

If $J_1\cap J^{\prime}_1=\emptyset$, then
$\sigma_{2,(j_3,j_4)}=\eta^{-1}\sigma_{2,(j_1,j_2)}\eta$, where
$\eta=(i_1,i_3)(i_2,i_4)\in A_d$.  We have $\eta^2=1$ and
$$\begin{array}{l}
(\sigma_{1,(j_1,j_2)}\sigma_{2,(j_3,j_4)})^r=
(\sigma_{1,(j_1,j_2)} \eta\sigma_{2,(j_1,j_2)}\eta^{-1})^r= \\
(\sigma_{1,(j_1,j_2)}
\sigma_{2,(j_1,j_2)}\eta^2)^r=(\sigma_{1,(j_1,j_2)}
\sigma_{2,(j_1,j_2)})^r=1.
\end{array}
$$
Therefore, the subgroup $H_2$ generated in $G$ by
$\sigma_{1,(j_1,j_2)}$ and $\sigma_{2,(j_3,j_4)}$ is isomorphic to
$D_r$. Note in addition that for each $i$ the group
$St_a(S_{d,i})=G\cap S_{d,i}$ is non-trivial and it is contains
either $\sigma_{i,(j_1,j_2)}$ or $\sigma_{i,(j_3,j_4)}$.

\subsection{Case $(2)$. Geometrical part}\label{case2part2}
We have $A_d= S_{d,i}\cap S_{d,j}$ for $i\neq j$. Denote by
$X_i=\widetilde X/S_{d-1,i}$, $\mathbb P^2=\widetilde X/S_{d,i}$,
and $X_0=\widetilde X/A_{d}$ the corresponding quotient spaces. They
can be arranged in the following commutative diagram (a fragment of
which is drown below) in which the morphisms $f_{0i}$, $i=1,\dots,
r$, are of degree two.

\begin{picture}(300,125)
\put(170,110){$\widetilde X$} \put(123,60){$X_i$}
\put(168,60){$X_0$} \put(211,60){$X_j$} \put(125,10){$\mathbb P^2$}
\put(213,10){$\mathbb P^2$}
\put(167,107){\vector(-1,-1){35}}\put(178,107){\vector(1,-1){35}}
\put(173,107){\vector(0,-1){33}}
\put(176,55){\vector(1,-1){33}}\put(172,55){\vector(-1,-1){33}}
\put(127,56){\vector(0,-1){32}} \put(215,56){\vector(0,-1){32}}
\put(154,85){$h_i$}\put(175,85){$h_0$}\put(202,85){$h_j$}
\put(116,38){$f_i$}\put(157,31){$f_{0i}$} \put(182,31){$f_{0j}$}
\put(217,38){$f_j$}
\end{picture}
\newline Since each $f_{0i}$ is conjugate to the covering $f_{01}:X_0\to \mathbb P^2$
branched along $B_1=B$, the covering $f_{0i}$ is branched over the
points of a cuspidal curve $B_i\subset \mathbb P^2$ having the same
degree and the same number of nodes and cusps as $B$ has. Denote by
$R_{i,0}\subset X_0$ the ramification curve of the covering
$f_{0i}$.

The group $D_r$ acts on $X_0$. The image of $\sigma_{i,(j_1,j_2)}\in
S_{d,i}$ (see subsection \ref{case2part1}) under the natural
epimorphism of $G$ to $D_r$ coincides with $\sigma_i$. Therefore,
the fixed point set of $\sigma_i$ coincides with $R_{i,0}$.

The surface $X_0$ is a normal projective variety. The set of its
singular points coincides with $f_{01}^{-1}(Sing B_1)$: over each
cusp we have a singular point of type $A_2$, and over each node we
have a singular point of type $A_1$. Therefore, for any $i$, $1\le
i\le r$, all the points of $f_{01}^{-1}(Sing B_1)$ belong to
$R_{i,0}$. In addition, as it follows from an observation made at
the end of subsection \ref{case2part1}, if two curves, say $R_{1,0}$
and $R_{2,0}$, meet at a nonsingular point $b\in X_0$, then the
point $b$ is common for all the curves $R_{i,0}$ and at this point
each pair of curves, $R_{i,0}$ and $R_{j,0}$, meet transversally.
Denote by $e$ the number of nonsingular points $b\in X_0$ common to
all the curves $R_{i,0}$.

Let $\nu:Z\to X_0$ be the minimal resolution of singularities. The
exceptional divisor of this resolution look as follows:
$$\overline E=\nu^{-1}(Sing X_0)=\bigcup_{
k=1}^c(\overline E_{1,s_k }\cup \overline E_{2,s_k })\cup \bigcup_{
l=2c+1}^{2c+n} \overline E_{s_ l},$$ where $\overline E_{1,s_k}$,
$\overline E_{2,s_k}$, $(1\le k\le c)$, are the irreducible
components of $\overline E$ contracted to the cusp $s_k$ of $X_0$,
and $\overline E_{s_l }$ ($ l=2c+1,\dots,2c+n$) is the irreducible
component of $\overline E$ contracted to the node $s_l$ of $X_0$.

Since $f_{0i}$ is a double covering branched along a cuspidal curve
$B_{i0}$, the above minimal resolution of singularities fits into
the following commutative diagram
$$
\xymatrix{ Z \ar[d]_{\overline f_{0i}} \ar[r]^{\nu } & X_0
\ar[d]^{f_{0i}} \\ \overline{\mathbb P}^2 \ar[r]_{\nu_i } & \mathbb
P^2, }
$$
where $\nu_i$ blows up once each of the singular points of $B_i$,
and $\overline f_{0i}$ is a two-sheeted covering of
$\overline{\mathbb P}^2$ branched over the strict transform
$\overline B_i\subset\overline{\mathbb P}^2$ of $B_i$.

Here is a more explicit description convenient for counting
intersection numbers.

Let $s$ be a cusp of $B_i$. Denote by $E\subset \overline{\mathbb
P}^2$ the exceptional curve $\nu_i^{-1}(s)=\mathbb P^1$ of $\nu_i$
lying over $s$. Then $\overline B_i\subset \overline{\mathbb P}^2$
meets $E$ at one point, it is non-singular at this point and has
there a simple tangency to $E$. The lift $\overline
R_{i,0}=f^{-1}_{0i}(\overline B_i)$ is the ramification curve of
$f_{0i}$, it is non-singular and coincides with the proper transform
of $R_{i,0}$. While $f^{-1}_{0i}(E)$ splits into $\overline
E_{1,s}\cup \overline E_{2,s}\subset Z$, a union of two smooth
curves intersecting transversally, so that
$$(\overline E_{1,s}^2)_{Z}=(\overline E_{2,s}^2)_{Z}=-2$$
and
$$(\overline E_{1,s},\overline E_{2,s})_{Z}=
(\overline E_{1,s},\overline R_{i,0})_{Z}=(\overline
E_{2,s},\overline R_{i,0})_{Z}=1.$$

Let $s$ be a node of $B_i$. Then $\overline B_i\subset
\overline{\mathbb P}^2$ meets the exceptional curve
$E=\nu_i^{-1}(s)$ at two points, $\overline B_i$ is nonsingular at
these points, and intersects $E$ transversally. The lift
$f^{-1}_{0i}(E)=\overline E_{s}\subset Z$ is the exceptional curve
of $\nu$, it meets $\overline R_{i,0}=f^{-1}_{0i}(\overline B_i)$,
which is non-singular, transversally, so that $(\overline
E_{s}^2)_{Z}=-2$ and $ (\overline E_{s},\overline R_{i,0})_{Z}=2$.

Now, let us show that $ (\overline R_{i,0},\overline R_{j,0})_{Z}$
does not depend on $i$ and $j$ if $i\neq j$. Consider the
commutative diagram
$$
\xymatrix{ \overline X \ar[d]_{\overline h_{0}} \ar[r]^{\overline
\mu } & \widetilde X \ar[d]^{h_{0}} \\ \overline Z \ar[r]_{\mu } &
X_0, }
$$
where $\overline X=\widetilde X\times_{X_0}\overline Z$ is the fibre
product of $\widetilde X$ and $\overline Z$ over $X_0$, and $\mu:
\overline Z\to X_0$ is the composition of $\nu$ and the blow ups of
the all intersection points of $(-2)$-curves $\overline E_{1,s_k}$
and $\overline E_{2,s_k}$ ($k=1,\dots,c$) lying in $Z$ (these curves
are those components of divisor $\overline E$ which are constructed
by $\nu$ to the cusps of $X_0$). Denote by $\overline
E_{1,2,s_k}\subset\overline Z$ the exceptional curve lying over the
point $\overline E_{1,s_k}\cap \overline E_{2,s_k}$ and, to simplify
the notation, denote by the same symbols, $\overline R_{i,0}$,
$\overline E_{1,s_k}$ and $\overline E_{2,s_k}$, and $\overline
E_{s_l}$, the strict transforms in $\overline Z$ of the curves
$\overline R_{i,0}$, $\overline E_{1,s_k}$ and $\overline E_{2,s_k}$
($k=1,\dots,c$), and $\overline E_{s_l}$ ($l=2c+1,\dots,2c+n$) in
$Z$.

It is easy to see that $\overline h_{0}$ is a Galois covering
branched along the curves $\overline E_{\cdot, s_k}$ and $\overline
E_{s_l}$. The ramification indices over the curves $\overline
E_{\cdot, s_k}$ are equal to $3$ (cf. local calculations in
\cite{Ku-Ku}, $\S 2$) and the ramification indices over the curves
$\overline E_{s_l
}$ are equal to $2$. The morphism $\overline \mu$ blows up once each
of the points lying over the nodes of $X_0$ and performs three blow
ups at each of the points lying over the cusps of $X_0$. Therefore,
the strict transforms $\overline \mu^{-1}(\widetilde
R_{i,(j_1,j_2)})$, $1\leq i\leq r$, $1\leq j_1,j_2\leq d$, pairwise
do not meet. But, $\bigcup_{j_1,j_2}\overline \mu^{-1}(\widetilde
R_{i,(j_1,j_2)})=\overline h_{0}^{-1}(\overline R_{i,0})$. Hence,
after the blow down of all curves $\overline E_{1,2,s_k}$, we get
$(\overline R_{i,0},\overline R_{j,0})_Z=c+e$ for $i\neq j$ and
these intersection numbers do not depend on $i$ and $j$. Note also
that the intersection numbers of the curves $\overline R_{j,0}$ and
an irreducible component of $\overline E$ do not depend on $j$.

The action of $D_r$  on $X_0$ lifts to an action on $Z$. The curve $\overline
R_{i,0}\subset Z$ (respectively, $R_{i,0}\subset X_0$) is the set of
fixed points of $\sigma_i\in D_r$. Since
$\sigma_i^{-1}\sigma_j\sigma_i\neq \sigma_j$ for $j\neq i$, the
curve $\sigma_i(\overline R_{j,0})\neq \overline R_{j,0}$
(respectively, $\sigma_i(R_{j,0})\neq R_{j,0}$) for $j\neq i$. In
particular, $R_{3,0}=\sigma_1(R_{2,0})\neq R_{2,0}$ and hence
$R_{2,0}+R_{3,0}=f_{01}^{-1}(D)$ for some curve $D\subset \mathbb
P^2$.

Since $D_r$ acts transitively on the set of curves $\overline
R_{i,0}$ (respectively, on the set of curves $R_{i,0}$), we have
$(\overline R_{1,0}^2)_Z=(\overline R_{2,0}^2)_Z=(\overline
R_{3,0}^2)_Z$ and it was shown that
$$(\overline R_{1,0},\overline R_{2,0})_Z=
(\overline R_{1,0},\overline R_{3,0})_Z= (\overline
R_{2,0},\overline R_{3,0})_Z.$$

Denote by $L$ the subspace of the Neron-Severi group $NS(Z)\otimes
\mathbb Q$ orthogonal (via the intersection form) to the subspace
$V_E$ generated by $\overline E_{1,s_k}$, $\overline E_{2,s_k}$,
$k=1,\dots,c$, and $\overline E_{s_l}$, $l=2c+1,\dots,2c+n$. The
intersection form is negative definite on $V_E$. Therefore, by Hodge
index theorem, the intersection form on $L$ has signature $(1,\dim
L-1)$.

In what follows we make, first, certain calculations  of
intersection numbers of some divisors in $L$. We project the
Neron-Severi classes of the divisors $\overline R_{i,0}$ to $L$,
denote the projections by $(\overline R_{i,0})_L$ and their
intersections in $L$ by $(\overline R_{i,0}\cdot \overline
R_{j,0})_L$ (the latter numbers, indeed, are equal to the
corresponding $\mathbb Q$-intersection numbers on the $\mathbb
Q$-variety $X_0$).

Observe, first, that $f^*_{i,0}B_i=2R_{i,0}$ and
$\nu^*R_{i,0}=2\overline R_{i,0}\mod L$. Thus,
$$
(\overline R_{i,0}^2)_L=
(\overline R_{j,0}^2)_L=\frac{1}{2}(\deg B)^2>0$$
for all $i,j$. Let $(\overline R_{2,0})_L=\lambda (\overline
R_{1,0})_L+T$, where $T\in L$ is orthogonal to $(\overline R_{1,0})_L$. We
have
$$
(\overline R_{2,0}\cdot \overline R_{1,0})_L=\lambda
(\overline R_{1,0}^2)_L=
(\overline R_{3,0}\cdot \overline R_{1,0})_L=
(\overline
R_{2,0}\cdot \overline R_{3,0})_L,$$ since the intersection numbers
of the curves $\overline R_{j,0}$ and an irreducible component of
$\overline E$ do not depend on $j$ and the intersection numbers
$(\overline R_{i,0},\overline R_{j,0})_Z$ also do not depend on
$i$ and $j$ for $i\neq j$.

Next, $\overline R_{2,0}+\overline R_{3,0}$ coincides with
$\nu^*(f^*_{01}(D))$. Therefore, $(\overline R_{2,0}+\overline
R_{3,0})_L$ is proportional to $(\overline
R_{1,0})_L=\frac{1}{2}\nu^*(f^*_{01}(B_1))$. Hence $(\overline
R_{3,0})_L=\lambda (\overline R_{1,0})_L-T$ and $\lambda >0$. We
have $(\overline R_{2,0}^2)_L=\lambda^2 (\overline R_{1,0}^2)_L+T^2=
(\overline R_{1,0}^2)_L$, therefore
$$T^2=(1-\lambda^2)(\overline R_{1,0}^2)_L\leq 0.$$
Hence $\lambda\geq 1$ and $\lambda= 1$ if and only if $T^2=0$,
that is, if and only if $T=0\in L$.

Since $(\overline R_{2,0}\cdot \overline R_{3,0})_L=(\overline
R_{2,0}\cdot \overline R_{1,0})_L$, we have
$$
\lambda^2  (\overline R_{1,0}^2)_L-T^2=\lambda (\overline R_{1,0}^2)_L
$$
and
$$
T^2=(\lambda^2-\lambda)(\overline R_{1,0}^2)_L\leq 0,$$ which yields
$\lambda\leq 1$. Combining this with the previous observations, we
get $\lambda = 1$ and $T=0$, which implies that $(\overline
R_{i,0})_L=(\overline R_{j,0})_L$ for any $i,j$. It allows us to
conclude that $\deg D=\deg B$, since $2\deg
D=((R_{2,0}+R_{3,0})^2)_L=4(R_{2,0}^2)_L=2\deg B$.

Therefore, we have $$\begin{array}{ll} (\overline R_{1,0}^2)_Z= &
(\frac{1}{2}\nu_1^*(B)-\sum(E_{1,s_k}+E_{2,s_k})-\sum
E_{s_l})^2=  \\
& \frac{1}{2} (\deg B)^2 -2c-2n=\frac12(\overline R_{1,0}, \overline
R_{2,0}+\overline R_{3,0})_Z\geq c . \end{array}$$

Hence
$$(\deg B)^2-4n\geq 6c.$$
On the other hand, as is known (see the proof of Lemma 3 in
\cite{Ku}),
$$
(\deg B)^2-2n\leq 6c
$$
for any generic covering of degree $d\ge 3$, while due to \cite{Ku1}
$n>0$ if $d>6$. The contradiction between these bounds eleminates
Case $(2)$.

\subsection{ Case $(3)$.} The symmetric group $G=S_{d+1}$ acts as the
permutation group on the set $I=\{ 1,\dots, d+1\}\subset \mathbb
N$. Denote by $H_i=\{ \gamma \in S_{d+1}\, |\, \gamma(i)=i\}$, so
that our $S_d=H_{d+1}$.

As in the proof of Theorem \ref{main1}, consider the quotient space
$\widetilde X/G=Y$ and the quotient map $\overline f: \widetilde
X\to Y$. The surface $Y$ is a normal projective variety. The
morphism $\overline f$ factors through $\widetilde f_i$, so that
$\overline f$ is the composition of the following morphisms

$$
\xymatrix{
\widetilde X  \ar[r]^{h } & X
\ar[r]^{f } & \mathbb P^2 \ar[r]^{r } & Y,  }
$$
where $r$ is a finite morphism of degree $d+1$. Since $S_{d}$ and
$S_{d+1}$ have no common normal subgroups, $\overline f$ is the
Galois expansion of $r$.

Let $\overline B\subset Y$ be the branch locus of $r$. We have
$r(B)=B_1\subset \overline B$. The preimage $r^{-1}(B_1)$ is the
union of  $B$ and some curve $B^{\prime}\subset \mathbb P^2$.

Since $Y$ is a normal projective surface, we can find a non-singular
projective curve $L\subset Y\setminus Sing(Y)$ which intersects
$\overline B$ transversally. Let $E=r^{-1}(L)$, $F=f^{-1}(E)$, and
$\widetilde F=\widetilde f^{-1}(E)$. Then $f_{\mid F}:F\to E$ is a
generic covering branched over $B\cap E$,  $\widetilde f_{\mid
\widetilde F}:\widetilde F\to E$ is the Galois expansion of the
generic covering $f_{\mid F}:F\to E$ with $Gal(\widetilde F/E)=S_d$,
and $\overline f_{\mid \widetilde F}:\widetilde F\to L$ is the
Galois expansion of the  covering $r_{\mid E}:E\to L$ with
$Gal(\widetilde F/L)=S_{d+1}$.

Consider the image $b_1=r(b)$ of a point $b\in B\cap E$. As in the
proof of Theorem \ref{main1}, it is easy to see that $r^{-1}(b_1)$
consists of $d-1$ points belonging to $B$ and one point point
belonging to $B^{\prime}$ (the ramification point of $r_{\mid E}$).
In other words, the covering $r_{\mid B^{\prime}}: B^{\prime}\to
B_1$ is of degree $1$, the covering $r_{\mid B}: B\to B_1$ is of
degree $d-1$, and $r^{*}(B_1)=B+ 2B^{\prime}$. In particular,
$$\begin{array}{ll}
\deg B^{\prime}\cdot\deg E= & ( B^{\prime},E)_{\mathbb P^2}=(B_1,L)_Y,  \\
\deg B\cdot\deg E= & ( B,E)_{\mathbb P^2}=(d-1)(B_1,L)_Y.
\end{array}
$$
Therefore
$$\deg B=(d-1)\deg B^{\prime}.$$
It follows, since $d=m^2K_X^2$ and $\deg B=m(3m+1)K_X^2$, that
$mK_X^2 +3$ is divisible by $m^2K_X^2-1$, which contradicts the
assumption $m^2K_X^2\geq  2\cdot 84^2$.

\subsection{ Case $(4)$}
Denote the standard imbedding $S_d\to  A_{d+2}$ by $\alpha$. For
each transposition $\sigma \in S_d$, the set $\widetilde
X_{\sigma}=\widetilde X_{\alpha(\sigma)}\subset \widetilde X$ of
fixed points of $\sigma$ is a nonsingular curve. Hence, for each
$\tau\in A_{d+2}$ which is conjugate to $\alpha(\sigma)$ the set
$\widetilde X_{\tau}$ of fixed points of $\tau$ is also a
nonsingular curve.

By \cite{Ku1}, if $d>6$, then the branch curve $B\subset \mathbb
P^2$ has at least one node. Therefore for each product
$\eta=\sigma_1\sigma_2$ of two commuting transpositions
$\sigma_1,\sigma_2\in S_d$, the set $\widetilde X_{\eta}$ of fixed
points of $\eta$ is finite and non-empty. It implies that for any
$\eta^{\prime}$, conjugated to $\alpha(\eta)$ in $ A_{d+2}$, the set
$\widetilde X_{\eta^{\prime}}$ is also finite and non-empty. On the
other hand, if $\sigma_1=(j_1,j_2)$ and $\sigma_2=(j_3,j_4)$, then
$\alpha(\sigma_1)=(j_1,j_2)((d+1),(d+2))$,
$\alpha(\sigma_2)=(j_3,j_4)((d+1),(d+2))$, and
$\alpha(\eta)=(j_1,j_2)(j_3,j_4)$ are conjugate to each other in  $
A_{d+2}$. Contradiction.

\section{Few applications}\label{applic}
\subsection{Deformation stability} The aim of this subsection is to prove a certain
deformation stability of examples given by Theorems \ref{main1} and
\ref{main2}. To state the corresponding results we need to fix few
notions. Namely, by a {\it $G$-manifold} we will mean a non-singular
projective manifold equipped with a regular action of the group $G$,
and by a smooth {\it $G$-family}, or {\it $G$-deformation}, of
$G$-manifolds we will mean a proper smooth morphism ({\it i.e.}, a
proper submersion) $p: \mathcal X \to B$, where $\mathcal X$ and $
B$ are smooth quasi-projective varieties, and $\mathcal X$ is
equipped with a regular action of $G$ preserving each fiber of $p$
(preservation of fibers means $p\circ G=p$).

\begin{prop}\label{stability}
If one of the fibers of a smooth $G$-family is the Galois expansion
of a generic covering of $\mathbb P^n$, then the whole family is
constituted from the Galois expansions of generic coverings of
$\mathbb P^n$.
\end{prop}

To prove Proposition \ref{stability} we need the following Lemma.
\begin{lem} \label{gener}  Let $p:\mathcal X \to
B$ be a smooth $G$-deformation of $G$-manifolds, $G$ being a finite
group. Then for each element $g\in G$ we have:
\begin{itemize} \item[$(i)$] the set $\mathcal
X^g=\{ x\in X\mid g(x)=x\}$ of fixed points of $g$ is a smooth
closed submanifold of  $\mathcal X$; \item[$(ii)$] the restriction
of $p$ to $\mathcal X^g$ is a smooth proper surjective morphism;
\item[$(iii)$] the intersection of $\mathcal X^g$ and each fibre
$X_t, t\in B$ of $p$ is transversal.
\end{itemize}
\end{lem}

\proof As is known, at any point $x\in \mathcal X$ fixed by a
subgroup $H$ of $G$, the action of $H$ can be linearized, which
means an existence of local analytic coordinates with respect to
which the action of $H$ is linear.

Let us resume Cartan's linearization procedure, see \cite{Cartan}.
Start from any system of local coordinates $z_1,\dots,z_n$ taking
value $0$ at a chosen point $x$ fixed by $H$. For any $h\in H$
denote by $h'$ the linear part of the Taylor expansion (with respect
to $z_1,\dots, z_n$) of $h$ at $x$. Then the change of coordinates
defined by the map $\sigma=\frac{1}{\vert H\vert}\sum_{g\in
H}(g')^{-1}g$ makes linear the action of $H$. Namely, it conjugates
$h$ and $h'$ for any $h\in H$, since $\sigma\circ h=h'\circ \sigma$
(indeed, $\sigma\circ h=\frac{1}{\vert H\vert}\sum_{g\in
H}(g')^{-1}g\circ h= \frac{1}{\vert H\vert}\sum_{g\in H}h'(g'\circ
h')^{-1}g\circ h= \frac{1}{\vert H\vert}h'\sum_{e\in
H}(e')^{-1}e=h'\circ\sigma$).

This change of coordinates is tangent to identity and it acts as
identity on each linear, with respect to $z_1,\dots,z_n$, subspace
on which $H$ acts already linear. Therefore, to prove $(i)$ it is
sufficient to linearize the action of $g$ (then, in new
coordinates the set $\mathcal X^g$ becomes linear), and to prove
$(ii)$ and $(iii)$ it is sufficient to pick any system of local
coordinates at $t\in B$ and include their lift into a system of
local coordinates $z_1,\dots,z_n$ (thus one gets for granted the
surjectivity of the projection at the level of tangent spaces, $
T_x(\mathcal X^g)\to T_{p(x)}B$; the properness and surjectivity
of $p: \mathcal X^g\to
B$ then follow from properness of $p:\mathcal X\to
B$ and closeness of $\mathcal X^g$ in $\mathcal X$). \qed

{\it Proof of Proposition {\rm \ref{stability}}.} We will give the
proof Proposition {\rm \ref{stability}} only in the case $n\leq 2$,
since the proof in the general case is similar.

Let $X_o, o\in B,$ be a fiber of $p$ which is the Galois expansion
of a generic covering $X_o \to\mathbb P^n$. The ramification locus
of this covering is a union of smooth codimension one manifolds $R_{
o (i,j)} $, $1\le i<j\le d$; these latter manifolds are the fixed
point sets of the transpositions $(i,j)\in S_d$. Due to Lemma
\ref{gener}, the fixed point sets $\mathcal X^ {(i,j)}\subset
\mathcal X$ of the same transpositions acting in $\mathcal X$ are
also smooth codimension
one manifolds, and for any $t\in
B$ the intersection $R_{ t (i,j)}= \mathcal{X}^ {(i,j)}\cap X_t$ is
transversal for each $(i,j)\in S_d$. Besides, if $\mathcal X^g\neq
\emptyset$ for some $g\in S_d$, $g\neq 1$, then, by Lemma
\ref{gener}, $\mathcal X^g\cap X_o\neq \emptyset$, and since the
action of $S_d$ on  $X_o$ is generic, it implies that $g$ is a
transposition in case $n=1$ and $g$ is either a transposition, or a
product of two disjoint transpositions, or a cyclic permutation of
length three in case  $n=2$.

 Let  $n=2$ and $g$ be a cyclic permutation
$(j_1,j_2,j_3)$ (the other cases can be treated in a similar way).
Let us show that the action of $S_{\{ j_1,j_2,j_3\}}$ on $\mathcal
X$, as well as on each of the fibers $X_t$, $t\in B$, is generic.
Indeed, without loss of generality we can assume that $\dim B=1$.
Then, by Lemma \ref{gener}, $\mathcal X^g$ is a smooth curve,
$\mathcal X^{(j_1,j_2)}$ and $\mathcal X^{(j_1,j_3)}$ are smooth
surfaces, and they all meet the fibers transversally. Since
$\mathcal X^g\cap X_o= R_{o(j_1,j_2)}\cap R_{o(j_1,j_3)}$,  by
applying once more Lemma \ref{gener}, we obtain that $\mathcal
X^{(j_1,j_2)}\cap \mathcal X^{(j_1,j_3)}= \mathcal X^g$ and
$X^g_t=R_{t(j_1,j_2)}\cap R_{t(j_1,j_3)}$ for any $t\in B$.
Therefore $\mathcal X^g$ (respectively, $X^g_t$) coincides with the
fixed point set under the action of $S_{\{ j_1,j_2,j_3\}}$ in
$\mathcal X$ (respectively, in $X_t$).

As a result, the action of $S_d$ on $\mathcal X$ and in each of
$X_t, t\in B$, is generic. Hence, the factor-space $\mathcal X/S_d$
is a smooth manifold and the induced morphism $p_1: \mathcal X/S_d\to
B$ is smooth and proper. Thus, it remains to notice that
$X_t/S_d=X_o/S_d=\mathbb P^n$ for any $t\in B$, since due to a
projective manifold $M$ is isomorphic to $\mathbb P^n$ as soon as
there exists a $C^\infty$-diffeomorphism $M\to\mathbb P^n$ which
maps the canonical class to the canonical class (for $n=1$ it is
known till Riemann; for $n=2$ one can use the Enriques-Kodaira
classification, see \cite{BPV}; it may be worth mentioning also
that, in fact, due to Siu \cite{Siu} in any dimension every compact
complex manifold deformation equivalent to $\mathbb P^n$ is
isomorphic to $\mathbb P^n$).\qed

\begin{cor}
$G$-varieties like in Theorems {\rm \ref{main1}} and {\rm
\ref{main2}} form connected components in the moduli space of,
respectively, $G$-curves and $G$-surfaces of general type. These
components are saturated {\rm (see the definition in Introduction)}.
In dimension $1$, $G$-varieties like in Theorem {\rm \ref{main1}}
also form proper subvarieties in the moduli space of curves of
general type.
\end{cor}

\begin{proof}
The first statement follows from Proposition \ref{stability}. The
second statement follows from the first one and Theorems {\rm
\ref{main1}} and {\rm \ref{main2}}. The third statement follows from
the first one and an observation that any birational transformation
of a one-parameter deformation family of genus $g\ge 2$ curves which
preserves each fiber and regular at the all points of all the
fibers, except a finite collection of fibers, extends indeed to a
transformation regular everywhere.
\end{proof}

\subsection{Examples of Dif$\ne $Def complex $G$-manifolds }
Here we consider regular actions of finite groups on complex
surfaces and construct diffeomorphic actions which are not
deformation equivalent. The idea is to pick diffeomorphic, but not
deformation equivalent, surfaces and apply to them Theorem
\ref{main2}.

Let $X$ be a rigid non real minimal
surface of general type, that is a minimal surface of general type
which is stable under deformations and not isomorphic to its own
conjugate, $\bar X$. Such surfaces are found in \cite{Kh-Ku2}.
Denote by $Y_1=\widetilde X$ the Galois expansion of a generic
$m$-canonical covering $X\to\mathbb P^2$, and by $Y_2$ its
conjugate, $Y_2=\bar Y_1$.

\begin{prop}\label{example1} Let $Y_1$ and $Y_2$ be as above
and let $m$ be like in Theorem {\rm \ref{main2}}. Then the actions
of $S_d=\Aut Y_1= \Aut Y_2$ on $Y_1$ and $Y_2$ are diffeomorphic,
but $Y_1$ and $Y_2$ are not $S_d$-deformation equivalent.
\end{prop}

\begin{proof} According to Theorem \ref{main2},
$\Aut Y_1= \Aut Y_2=S_d$ where $d$ is the degree of
the $m$-canonical covering $X\to\mathbb P^2$.
The action of $S_d$ in $Y_1$ is tautologically diffeomorphic to
that in $Y_2$, since $Y_2=\bar Y_1$.

Assume that $Y_1=\widetilde X$ and $Y_2=\bar Y_1$ are
$S_d$-deformation equivalent. Let $p:\mathcal X \to
B$ be a smooth $S_d$-deformation connecting them (the treatment of a
chain of deformation families is literally the same). By Proposition
\ref{stability}, for any $t\in B$ the covering $X_t\to\mathbb P^n$
is generic. Hence, $\mathcal X/S_{d-1}\to B$ is a deformation family
connecting $X= Y_1/S_{d-1}$ with $\bar X=Y_2/S_{d-1}$, which is a
contradiction.
\end{proof}

\subsection{Examples of Dif$\ne $Def real $G$-manifolds} Here we extend the
category of  $G$-manifolds, namely, we consider finite subgroups of
the Klein extension of the automorphism group. Let us recall that
the Klein group ${\operatorname Kl} (X)$ of a complex variety $X$
is, by definition, the group consisting of biregular isomorphisms
$X\to X$ and $X\to\bar X$ (some people call it the group of
dianalytic automorphisms). If $X$ is a real manifold and $c$ its
real structure, then there is an exact sequence $$1\to \langle
c\rangle={\mathbb Z}/2\to \Kl(X)\to \Aut X\to 1.$$

Pick two real Campedelli surfaces $(X_1,c_1)$, $(X_2,c_2)$
constructed in \cite{Kh-Ku1}, Section 2. As is shown in
\cite{Kh-Ku1}, these particular surfaces are not real deformation
equivalent, but their real structures, $c_1:X_1\to X_1$ and
$c_2:X_2\to X_2$, are diffeomorphic.

The Campedelli surfaces are minimal surfaces of general type. Thus,
we can consider $m$-canonical generic coverings $X_1\to {\mathbb
P^2}$ and $X_2\to{\mathbb P^2}$. Moreover, we can choose these
coverings to be real, that is to be equivariant with respect to the
usual, complex conjugation, real structure on ${\mathbb P^2}$ and
the real structures $c_1, c_2$ on $X_1,X_2$. Denote by $\widetilde
X_1\to {\mathbb P^2}$ and $\widetilde X_2\to {\mathbb P^2}$ the
Galois expansions. The surfaces
$\widetilde X_1$ and $\widetilde X_2$ are real with real structures
lifted from $\mathbb P^2$.

\begin{prop} The Klein groups $\Kl(\widetilde X_1)$ and $\Kl(\widetilde X_2)$
are isomorphic, their
actions are diffeomorphic, while there exists no equivariant
deformation connecting $(\widetilde X_1,\Kl(\widetilde X_1))$ with
$(\widetilde X_2, \Kl(\widetilde X_2))$.
\end{prop}

\begin{proof} As in the proof of Proposition \ref{example1},
the non existence of an equivariant deformation is a straightforward
consequence of Proposition \ref{stability}.

Due to Theorem \ref{main2} to prove the other two statements, it is
sufficient to construct a real diffeomorphism $\widetilde
X_1\to\widetilde X_2$ respecting the Galois action. In its turn, to
reach this task, it is sufficient to construct real diffeomorphisms
$X_1\to X_2$, $\mathbb P^2\to\mathbb P^2$ commuting with the initial
$m$-canonical generic coverings, $X_1\to \mathbb P^2$ and $X_2\to
\mathbb P^2$.

Now, we need to recall some details of the construction of surfaces
$X_1$, $X_2$ from \cite{Kh-Ku1}. The construction starts from a real
one-parameter family of Campedelli line arrangements $\mathcal
L(t)$, $t\in T =\{|t|\le 1, t\in\mathbb C\} $, consisting of seven
lines $L_1(t),\dots L_7(t)$ labeled by non-zero elements $\alpha\in
(\mathbb Z/2\mathbb Z)^3$. The lines are real for real values of
$t$, and the family performs a triangular transformation at $t=0$
(see the definition of triangular transformation in  \cite{Kh-Ku1}).
Consider the Galois covering $Y\to \mathbb P^2\times T$ with Galois
group $(\mathbb Z/2\mathbb Z)^3$ branched in $\sum^7_1
\mathcal L_i, \mathcal L_i=\{(p,t)\in \mathbb P^2\times T
: p\in L_i(t)\}$ and defined by the chosen labelling of the lines.
The fibers $Y_t$ under the projection of $Y$ to $T$ are nonsingular
Campedelli surfaces for generic $t$, in particular, for any $t\ne 0$
close to $0$. The surfaces $X_1, X_2$ we are interested in are given
by $Y_t$ with, respectively, positive and negative $t$ close to $0$.
The fiber $Y_0$ has two singular points; each of these points is a
so-called $T(-4)$-singularity; these points are not real but complex
conjugate to each other. For each nonsingular Campedelli surface
$Y_t$, $t\ne 0$, and each $1\le i\le 7$, the pull-back
$L^*_i(t)\subset Y_t$ of $L_i(t)\subset \mathbb P^2$ represents the
bi-canonical class, $ [L^*_i(t)]=2K_{Y_t}$.

Since $Y\to \mathbb P^2\times T$ is a finite morphism, the divisors
$E_m= m[\mathcal L_i^*]$ are relatively very ample for any $m$
sufficiently large (and any $i$), see for example
\cite{HartshorneLN}. Pick such an integer $m$ and consider the real,
relative to $T $, imbedding of $Y$ into $\mathbb P^N\times T$
defined by the linear system $|E_m|$ (to show existence of relative
to $T $ global sections, one can twist $E_m$ by a pull-back of a
very ample divisor on $T $). Note, that according to
$[L^*_i(t)]=2K_{Y_t}$ it defines $2m$-canonical imbedding of
Campedelli surfaces $Y_t$ to $\mathbb P^N$. As it follows from
Theorem 0.1 in \cite{Ku-Ku},  if $m\geq 5$, then the projection
$Y_t\to\mathbb P^2$ from a generic $\mathbb P^{N-3}$ is a generic
covering for any but finite number of $t\in T $, in particular, for
any $t\ne 0$ close to $0$. For a real $\mathbb P^{N-3}$ the
projection $Y\to \mathbb P^2\times T $ is real and, for a
sufficiently generic real $\mathbb P^{N-3}$, the two singular points
of $Y_0$ project to two distinct complex conjugate points. Denote
the singular points by $y$, $\bar y$, and their projections by $b$,
$\bar b$. One can show, repeating word-by-word the arguments from
\cite{Ku-Ku}, that the projection $Y_0\to \mathbb P^2$ is generic
everywhere (generic at singular point means that the fibre of
projection passing through the singular point $y\in Y_0$ (resp.
$\bar y\in Y_0$) is in generic position with respect to the tangent
cone $C_bY_0$ (resp. $C_{\bar b}Y_0$)).

Restrict, now, our attention to small values of $t$. The coverings
$Y_t\to\mathbb P^2$ are generic for $t\ne 0$, and the branching
curves $B_t$ are cuspidal. For $t=0$, the branching curve $B_0$ is
cuspidal everywhere, except two distinct, complex conjugate, points,
$b$ and $\bar b$. Cut out small, Milnor, complex conjugate balls
$V(b), V(\bar b)$ around these two points. Use a family of
Morse-Lefschetz diffeomorphisms to complete the isotopy
$B_{te^{i\phi}}\setminus (V(b) \cup V(\bar b))$ by an isotopy inside
$V(b)$, and then complete it by a complex conjugate isotopy inside
$V(\bar b)$. This isotopy provides an equivariant diffeomorphism
between Galois coverings branched in $B_t$ and, respectively,
$B_{-t}$. The Morse-Lefschetz diffeomorphisms can be seen as
transformations acting identically on the complement of $V(b)\cup
V(\bar b)$. Therefore, the epimorphism $\pi_1(\mathbb P^2\setminus
B_t)\to S_d$ defining the Galois coverings is not changing, so that
the constructed diffeomorphism between the Galois coverings acts
from $\widetilde X_1$ to $\widetilde X_2$ and it is equivariant with
respect to the Galois action. It is also equivariant with respect to
the real structure. Thus, it remains to notice that due to Theorem
\ref{main2} the full automorphism groups $\Aut(\widetilde X_1)$ and
$\Aut(\widetilde X_2)$
coincide with the Galois group.
\end{proof}

\end{document}